\nonstopmode
\documentclass[]{amsart}
\usepackage{amssymb}
\usepackage{graphicx}
\usepackage{epsfig}
\allowdisplaybreaks

\def\undersetbrace#1\to#2{\underbrace{#2}_{#1}}
\def\oversetbrace#1\to#2{\overbrace{#2}^{#1}}
\def\AMSunderset#1\to#2{\underset{#1}{#2}}
\def\AMSoverset#1\to#2{\overset{#1}{#2}}

\def\East#1#2{\overset{#1}{\longrightarrow}}
\swapnumbers

\newtheorem{prop}[subsection]{Proposition}
\newtheorem*{prop*}{Proposition}
\newtheorem{thm}[subsection]{Theorem}
\newtheorem*{thm*}{Theorem}

\newtheorem*{lem*}{Lemma}

\newtheorem*{kor*}{Corollary}
\newenvironment{demo}[1]{\par\smallskip\noindent{\bf #1.}}{\par\smallskip}

\def\cit#1#2{\ifx#1!\cite{#2}\else#2\fi} 
\def\idx{}               
\def\ign#1{}             
\def\o{\circ}
\def\X{\mathfrak X}
\def\al{\alpha}
\def\be{\beta}
\def\ga{\gamma}

\def\ep{\varepsilon}

\def\et{\eta}
\def\th{\theta}

\def\la{\lambda}
\def\rh{\rho}

\def\ph{\varphi}

\def\ps{\psi}

\def\Ga{\Gamma}

\def\Om{\Omega}
\def\i{^{-1}}
\def\x{\times}
\def\p{\partial}
\let\on=\operatorname
\def\D{T} 
\def\L{\mathcal L}
\def\today{\ifcase\month\or
 January\or February\or March\or April\or May\or June\or
 July\or August\or September\or October\or November\or December\fi
 \space\number\day, \number\year}
\def\AMSonly#1{}
\begin{document}
\title[Vanishing geodesic distance]
{Vanishing geodesic distance on spaces of submanifolds
and diffeomorphisms
}
\author{Peter W. Michor, David Mumford}
\address{
Peter W. Michor:
Fakult\"at f\"ur Mathematik, Universit\"at Wien,
Nordbergstrasse 15, A-1090 Wien, Austria; {\it and:}
Erwin Schr\"odinger Institut f\"ur Mathematische Physik,
Boltzmanngasse 9, A-1090 Wien, Austria
}
\email{Peter.Michor@univie.ac.at}
\address{
David Mumford:
Division of Applied Mathematics, Brown University,
Box F, Providence, RI 02912, USA}
\email{David\_{}Mumford@brown.edu}
\thanks{Both authors were supported by NSF Grant 007~4276.
     PWM was supported by FWF Projects P~14195 and P~17108 and by Centre
     Bernoulli, Lausanne}

\begin{abstract} 
The $L^2$-metric or Fubini-Study metric on the non-linear
Grassmannian of all submanifolds of type $M$ in a Riemannian
manifold $(N,g)$ induces geodesic distance 0. 
We discuss another metric which involves the mean curvature and
shows that its geodesic distance is a good topological metric.
The vanishing phenomenon for the geodesic distance holds also
for all diffeomorphism groups for the $L^2$-metric. 
\end{abstract}
\def\LaTeXonly{}

\maketitle

\section{Introduction} \label{nmb:1}

In \cit!{M98} we studied the $L^2$-Riemannian metric on the space of all 
immersions $S^1\to \mathbb R^2$. This metric is invariant under the group 
$\on{Diff}(S^1)$ and we found that it induces {\it vanishing geodesic 
distance} on the quotient space $\on{Imm}(S^1,\mathbb R^2)/\on{Diff}(S^1)$.
In this paper we extend this result to the general situation
$\on{Imm}(M,N)/\on{Diff}(M)$ for any compact manifold $M$ and
Riemannian manifold $(N,g)$ with $\dim N>\dim M$. On the open subset 
$\on{Emb}(M,N)/\on{Diff}(M)$, which may be identified with the space
of all submanifolds of diffeomorphism type $M$ in $N$ (the non-linear 
Grassmanian or differentiable `Chow' variety) this says that the infinite 
dimensional analog of the Fubini Study metric induces vanishing geodesic 
distance. The picture that emerges for these infinite-dimensional manifolds 
is quite interesting: there are simple expressions for the Christoffel 
symbols and curvature tensor, the geodesic equations are simple and of 
hyperbolic type and, at least in the case of plane curves, the geodesic 
spray exists locally. {\it But} the curvature is positive and unbounded 
in some high frequency directions, so these spaces wrap up on themselves 
arbitrarily tightly, allowing the infimum of path lengths between two given
points to be zero. 

We also carry over to the general case the stronger metric from
\cite{M98} which weights the $L^2$ metric using the second fundamental 
form. It turns out that we have only to use the mean curvature in order 
to get positive geodesic distances, hence a good topological metric on 
the space $\on{Emb}(S^1,\mathbb R^2)/\on{Diff}(S^1)$. The reason is that the first variation of the volume of a submanifold depends on the mean curvature
and the key step is showing that the square root of the volume of $M$
is Lipschitz in our stronger metric. The formula for this metric is:
\begin{equation*}
G^{A}_f(h,k) :=
\int_{M}(1+A\|\on{Tr}^{f^*g}(S^f)\|_{g^{N(f)}}^2)g(h,k)\on{vol}(f^*g)
\end{equation*}
where $S^f$ is the second fundamental form of the immersion
$f$; section \ref{nmb:3} contains the relevant estimates. 
In section \ref{nmb:4} we also compute the sectional curvature
of the $L^2$-metric in the hope to relate the vanishing of the
geodesic distance to unbounded positivity of the sectional
curvature: by going through ever more positively curved parts
of the space we can find ever shorter curves between any two
submanifolds.

In the final section \ref{nmb:5} we show that the vanishing of
the geodesic distance also occurs on the Lie group of all
diffeomorphisms on each connected Riemannian manifold. Short
paths between any 2 diffeomorphisms are constructed by using
rapidly moving compression waves in which individual points are 
trapped for relatively long times. We compute the sectional curvature 
also in this case.

\section{The manifold of immersions} \label{nmb:2}

\subsection{Conventions} \label{nmb:2.1}

Let $M$ be a compact smooth connected manifold of dimension $m\ge 1$ and let
$(N,g)$ be a connected Riemannian manifold of dimension $n>m$. We
shall use the following spaces and manifolds of smooth
mappings.

\begingroup
\parindent=0cm
\everypar{\hangindent=1.5cm \hangafter=1} 

$\on{Diff}(M)$, the regular Lie group (\cit!{4},~38.4) of all
diffeomorphisms of $M$. 

$\on{Diff}_{x_0}(M)$, the subgroup of diffeomorphisms fixing
$x_0\in M$.

$\on{Emb}=\on{Emb}(M,N)$, 
the manifold of all smooth embeddings $M\to N$. 

$\on{Imm}=\on{Imm}(M,N)$, 
the manifold of all smooth immersions $M\to N$.
For an immersion $f$ the tangent space with foot point $f$
is given by 
$T_f\on{Imm}(M,N)=C^\infty_f(M,TN)=\Ga(f^*TN)$, the space of
all vector fields along $f$.

$\on{Imm}_f=\on{Imm}_f (M,N)$, 
the manifold of all smooth 
{\it free} immersions $M\to N$, i.e., those with trivial isotropy group
for the right action of $\on{Diff}(M)$ on $\on{Imm}(M,N)$.

$B_e=B_e(M,N)=\on{Emb}(M,N)/\on{Diff}(M)$, the manifold of
submanifolds of type $M$ in $N$, the base of a smooth principal
bundle, see \ref{nmb:2.2}.

$B_i=B_{i}(M,N)=\on{Imm}(M,N)/\on{Diff}(M)$, an infinite
dimensional `orbifold', whose points are, roughly speaking,
smooth immersed submanifolds of type $M$ in $N$,
see \ref{nmb:2.4}.

$B_{i,f}=B_{i,f}(M,N)=\on{Imm}_{f}(M,\mathbb
R^2)/\on{Diff}(M)$, a manifold, the base of a principal fiber bundle, see
\ref{nmb:2.3}.

\endgroup

For a smooth curve $f:\mathbb R\to C^\infty(M,N)$ corresponding
to a mapping
$f:\mathbb R\x M\to N$, we shall denote by $\D f$ the curve of
tangent mappings, so that $\D f(t)(X_x)=T_x(f(t,\quad)).X_x$.
The time derivative will be denoted by either 
$\p_t f=f_t:\mathbb R\x M\to TN$. 

\subsection{The principal bundle of embeddings 
$\on{Emb}(M,N)$} \label{nmb:2.2}

We recall some basic results whose proof can be found in {\cit!{4}):

{\it \noindent (A) The set $\on{Emb}(M,N)$ of all smooth embeddings 
$M\to N$ is an open subset of the smooth Fr\'echet manifold 
$C^\infty(M,N)$ of all smooth mappings $M\to N$ 
with the $C^\infty$-topology. 
It is the total space of a smooth principal bundle
$\pi:\on{Emb}(M,N)\to B_e(M,N)$ with structure group 
$\on{Diff}(M)$, the smooth regular Lie group group of all 
diffeomorphisms of $M$, whose base $B_e(M,N)$ is the smooth 
Fr\'echet manifold of all submanifolds of $N$ of type $M$, 
i.e., the smooth manifold of all simple closed curves in $N$.
{\rm (\cit!{4},~44.1)} 

\noindent (B) This principal bundle admits a smooth principal connection
described by the horizontal bundle whose fiber $\mathcal N_c$ over $c$ 
consists of all vector fields $h$ along $f$ such that
$g(h, Tf)=0$. The parallel transport for this connection
exists and is smooth. {\rm (\cit!{4},~39.1 and 43.1)}}

\subsection{Free immersions} \label{nmb:2.3}

The manifold $\on{Imm}(M,N)$ of all immersions $M\to N$
is an open set in the manifold $C^\infty(M,N)$ and thus itself a
smooth manifold. 
An immersion $f:M\to N$ is called {\it free} if $\on{Diff}(M)$
acts freely on it, i.e., $f\o \ph=c$ for $\ph\in \on{Diff}(M)$ implies
$\ph=\on{Id}$. 
We have the following results:
\begin{itemize}
\item
{\it If $\ph\in \on{Diff}(M)$ has a fixed point and if $f\o\ph=f$ for some 
immersion $f$ then $\ph=\on{Id}$.} This is (\cit!{1},~1.3).
\item
{\it If for $f\in \on{Imm}(M,N)$ there is a point $x\in f(M)$
with only one preimage then $f$ is a free immersion.} This is
(\cit!{1},~1.4). There exist free immersions without such
points.
\item
{\bf The manifold $B_{i,f}(M,N)$}
{\rm (\cit!{1},~1.5)}
{\it The set $\on{Imm}_{f}(M,N)$ of all free immersions
is open in $C^\infty(M,N)$ and thus a smooth submanifold. 
The projection 
\begin{displaymath}
\pi:\on{Imm}_{f}(M,N)\to 
\frac{\on{Imm}_{f}(M,N)}{\on{Diff}(M)}=:B_{i,f}(M,N)
\end{displaymath}
onto a Hausdorff smooth manifold is a smooth principal fibration with structure
group $\on{Diff}(M)$. 
By {\rm (\cit!{4},~39.1 and 43.1)} this fibration 
admits a smooth principal connection described by the
horizontal bundle with fiber $\mathcal N_c$ consisting of all vector fields $h$ 
along $f$ such that $g(h,Tf)=0$. This connection admits a
smooth parallel transport over each smooth curve in the base manifold.
}
\end{itemize}

We might view $\on{Imm}_{f}(M,N)$ as the nonlinear Stiefel 
manifold of parametrized submanifolds of type $M$ in $N$ and consequently 
$B_{i,f}(M,N)$ as the nonlinear Grassmannian of unparametrized
submanifolds of type $M$.

\subsection{Non free immersions} \label{nmb:2.4}

Any immersion is proper since $M$ is compact and thus 
by (\cit!{1},~2.1) the orbit space 
$B_i(M,N)=\on{Imm}(M,N)/\on{Diff}(M)$ is Hausdorff. 
Moreover, by (\cit!{1},~3.1 and 3.2) for any immersion $f$ 
the isotropy group $\on{Diff}(M)_f$ 
is a finite group which acts as group of covering transformations for a
finite covering $q_c:M\to \bar M$ such that $f$ factors over $q_c$ to a free
immersion $\bar f:\bar M\to N$ with $\bar f\o q_c = f$. 
Thus the subgroup $\on{Diff}_{x_0}(M)$ of all diffeomorphisms $\ph$ fixing
$x_0\in M$ acts freely on $\on{Imm}(M,N)$.
Moreover, for each $f\in \on{Imm}$  the submanifold
$\mathcal Q(f)$ from \ref{nmb:4.4}, \thetag{1}
is a slice in a strong sense:
\begin{enumerate}
\item[$\bullet$] 
   $\mathcal Q(f)$ is invariant under the isotropy group
   $\on{Diff}(M)_f$.
\item[$\bullet$] 
   If $\mathcal Q(f)\o \ph \cap \mathcal Q(f)\ne\emptyset$ for
   $\ph\in\on{Diff}(M)$ then $\ph$ is already in the isotropy group 
   $\ph\in\on{Diff}(M)_f$.
\item[$\bullet$] 
   $\mathcal Q(f)\o \on{Diff}(M)$ is an invariant open neigbourhood of
   the orbit $f\o\on{Diff}(M)$ in $\on{Imm}(M,N)$ which
   admits a smooth retraction $r$ onto the orbit. The fiber
   $r\i(f\o\ph)$
   equals $\mathcal Q(f\o\ph)$.
\end{enumerate}
Note that also the action 
$$
\on{Imm}(M,N)\x \on{Diff}(M) \to 
\on{Imm}(M,N)\x\on{Imm}(M,N),\qquad (f,\ph)\mapsto
(f,f\o\ph)
$$
is proper so that all assumptions and conclusions of Palais' slice theorem
\cit!{P} hold.
This results show that the orbit space 
$B_i(M,N)$ has only singularities of orbifold type 
times a Fr\'echet space.
We may call the space $B_i(M,N)$ an infinite dimensional {\it
orbifold}. The projection 
$\pi:\on{Imm}(M,N)\to B_i(M,N)
=\on{Imm}(M,N)/\on{Diff}(M)$
is a submersion off the singular points and has only mild singularities at the
singular strata.
The normal bundle $\mathcal N_f$ mentioned in \ref{nmb:2.2} is well defined
and is a smooth vector subbundle of the tangent bundle. 
We do not have a principal bundle and thus no principal connections, but we
can prove the main consequence, the existence of horizontal paths,
directly:

\begin{prop}\label{nmb:2.5}
For any smooth path $f$ in $\on{Imm}(M,N)$ there exists a smooth
path $\ph$ in $\on{Diff}(M)$ with $\ph(t,\quad)=\on{Id}_{M}$ depending
smoothly on $f$ such that
the path $h$ given by $h(t,\th)=f(t,\ph(t,\th))$ is horizontal:
$g(h_t, \D h)=0$.
\end{prop}

\begin{demo}{Proof}
Let us write $h=f\o \ph$ for $h(t,x)=f(t,\ph(t,x))$, etc. 
We look for $\ph$
as the integral curve of a time dependent vector field $\xi(t,x)$ on
$M$, given by $\p_t\ph=\xi\o \ph$.
We want the following expression to vanish:
\begin{align*}
g\bigl( \p_t(f\o\ph),\D(f\o\ph)\bigr)
&= g\bigl((\p_tf\o\ph + (\D f\o\ph).\p_t\ph,(\D f\o\ph).\D\ph\bigr)
\\&
=\bigl(g(\p_tf,\D f)\o\ph\bigr).\D\ph 
  +g\bigl( (\D f\o \ph)(\xi\o\ph),(\D f\o\ph).\D\ph\bigr)
\\&
=\bigl(\bigl(g(\p_tf,\D f) +g(\D f.\xi,\D f)\bigr)\o\ph\bigr).\D\ph
\end{align*}
Since $\D \ph$ is everywhere invertible we get
$$
0=g\bigl( \p_t(f\o\ph),\D(f\o\ph)\bigr) \iff 
0=g(\p_tf,\D f) +g(\D f.\xi,\D f)
$$
and the latter equation determines the non-autonomous vector
field $\xi$ uniquely. 
\qed\end{demo}

\subsection{Curvatures of an immersion}\label{nmb:2.6}
Consider a fixed immersion $f\in\on{Imm}(M,N)$. 
The {\it normal bundle} $N(f)=Tf^\bot\subset f^*TN\to M$ has
fibers $N(f)_x=\{Y\in T_{f(x)}N:g(Y,T_xf.X)=0\text{ for all
}X\in T_xM\}$.
Every vector field $h:M\to TN$ along $f$ then splits as
$h=Tf.h^\top+h^\bot$ into its tangential component
$h^\top\in\X(M)$ and its normal component $h^\bot\in
\Ga(N(f))$.

Let $\nabla^g$ be the Levi-Civita covariant derivative of $g$
on $N$ and let $\nabla^{f^*g}$ the Levi-Civita covariant
derivative of the pullback metric $f^*g$ on $M$. 
The \idx{\it shape operator} or \idx{\it second fundamental
form} $S^f\in \Ga(S^2 T^*M\otimes N(f))$ of $f$ is then given by
\begin{equation}
S^f(X,Y)=\nabla^g_X(Tf.Y)-Tf.\nabla^{f^*g}_XY\quad \text{  for }X,Y\in
\X(M).
\tag{1}\end{equation}
It splits into the following irreducible components
under the action of the group $O(T_xM)\x O(N(f)_x)$: the mean
curvature $\on{Tr}^{f^*g}(S^f)=\on{Tr}((f^*g)\i\o S^f)\in
\Ga(N(f))$ and the trace free shape operator 
$S^f_0 = S^f-\on{Tr}^{f^*g}(S^f)$. 
For $X\in \X(M)$ and $\xi\in\Ga(N(f))$, i.e., a normal vector field
along $f$, we may also split $\nabla^g_X\xi$ into the components
which are tangential and normal to $Tf.TM$,
\begin{equation}
\nabla^g_X\xi = -Tf.L^f_{\xi}(X) + \nabla^{N(f)}_X\xi 
\tag{2}\end{equation}
where $\nabla^{N(f)}$ is the induced connection in the normal
bundle respecting the metric $g^{N(f)}$ induced by $g$, 
and where the {\it Weingarten} tensor field
$L^f\in\Ga(N(f)^*\otimes T^*M\otimes TM)$ corresponds to the
shape operator via the formula 
\begin{equation}
(f^*g)(L^f_\xi(X),Y) = g^{N(f)}(S^f(X,Y),\xi).
\tag{3}\end{equation}
Let us also split the Riemann curvature $R^g$ into tangential
and normal parts: For $X_i\in \X(M)$ or $T_xM$ we have
({\it theorema egregium}):
\begin{align}
g(R^g&(Tf.X_1,Tf.X_2)(Tf.X_3),Tf.X_4) = (f^*g)(R^{f^*g}(X_1,X_2)X_3,X_4)
+\notag\\
&+ g^{N(f)}(S^f(X_1,X_3),S^f(X_2,X_4)) - g^{N(f)}(S^f(X_2,X_3),S^f(X_1,X_4)).
\tag{4}\end{align}
The normal part of $R^g$ is then 
       given by ({\it Codazzi-Mainardi equation}):
\begin{align}
(R^g&(Tf.X_1,Tf.X_2)(Tf.X_3))^\bot = 
\notag\\&
= \bigl(\nabla^{N(f)\otimes T^*M\otimes T^*M}_{X_1} S^f\bigr)(X_2,X_3) -
\bigl(\nabla^{N(f)\otimes T^*M\otimes T^*M}_{X_2} S^f\bigr)(X_1,X_3). 
\tag{5}\end{align}

\subsection{Volumes of an immersion}\label{nmb:2.7}
For an immersion $f\in\on{Imm}(M,N)$, we consider the volume density
$\on{vol}^g(f)=\on{vol}(f^*g)\in\on{Vol}(M)$ on $M$
given by the local formula
$\on{vol}^g(f)|_U=\sqrt{\det((f^*g)_{ij})}|du^1\wedge \dots \wedge du^m|$
for any chart $(U,u:U\to \mathbb R^m)$ of $M$, 
and
the induced volume function $\on{Vol}^g:\on{Imm}(M,N)\to \mathbb
R_{>0}$ which is given by $\on{Vol}^g(f)=\int_M\on{vol}(f^*g)$.
The tangent mapping of $\on{vol}:\Ga(S^2_{>0}T^*M) \to \on{Vol}(M)$
is given by 
$d\on{vol}(\ga)(\et)=\frac12 \on{Tr}(\ga\i.\et)\on{vol}(\ga)$.
We consider the pullback mapping 
$P_g:f\mapsto f^*g, P_g:\on{Imm}(M,N)\to
\Ga(S^2_{>0}T^*M)$. 
A version of the following lemma is \cit!{Ka},~1.6. 

\begin{lem*}
The derivative of 
$\on{vol}^g=\on{vol}\o P_g:\on{Imm}(M,N)\to \on{Vol}(M)$ is 
\begin{align*}
d\on{vol}^g&(h)=d(\on{vol}\o P_g)(h) 
\\&
=-\on{Tr}^{f^*g}(g(S^f,h^\bot))\on{vol}(f^*g)
+\tfrac12 \on{Tr}^{f^*g}(\mathcal
L_{h^\top}(f^*g))\on{vol}(f^*g).
\\&
=-(g(\on{Tr}^{f^*g}(S^f),h^\bot))\on{vol}(f^*g)
+\on{div}^{f^g}(h^\top)\on{vol}(f^*g).
\end{align*}
\end{lem*}

\begin{demo}{Proof}
We consider a curve $t\mapsto f(t,\quad)$ in $\on{Imm}$ with
$\p_t|_0f=h$.
We also use a chart  $(U,u:U\to \mathbb R^m)$ on $M$. 
Then we have 
\begin{align*}
f^*g|_U&=\sum_{i,j}(f^*g)_{ij}du^i\otimes du^j =
\sum_{i,j}g(Tf.\p_{u_i},Tf.\p_{u^j})du^i\otimes du^j
\\
\p_t \on{vol}^g(f)|_U 
&= \frac{\det((f^*g)_{ij})(f^*g)^{kl} \p_t(f^*g)_{lk}}
        {2\sqrt{\det((f^*g)_{ij})}}
   |du^1\wedge \dots \wedge du^m|
\end{align*}
where
\begin{align*}
\p_t(f^*g)_{ij}&=\p_t g(Tf.\p_{u^i},Tf.\p_{u^j})
\\&
= g(\nabla^g_{\p_t}(Tf.\p_{u^i}),Tf.\p_{u^j})
+g(\p_{u^i},\nabla^g_{\p_t}(Tf.\p_{u^j}))
\\
g(\nabla^g_{\p_t}Tf.\p_{u^i},Tf.\p_{u^j})
&= 
g(\nabla^g_{\p_{u^i}}Tf.\p_{t}+Tf.\on{Tor}+Tf.[\p_t,\p_{u^i}],Tf.\p_{u^j})
\\&
= g(\nabla^g_{\p_{u^i}}(Tf.\p_{t})^\bot,Tf.\p_{u^j})
+ g(\nabla^g_{\p_{u^i}}(Tf.\p_{t})^\top,Tf.\p_{u^j})
\\
g(\nabla^g_{\p_{u^i}}(\p_{t}f)^\bot,Tf.\p_{u^j})
&
=g(-Tf.L_{(\p_{t}f)^\bot}\p_{u^i},Tf.\p_{u^j})
+g(\nabla^{N(f)}_{\p_{u^i}}(\p_{t}f)^\bot,Tf.\p_{u^j})
\\&
=-(f^*g)(L_{(\p_{t}f)^\bot}\p_{u^i},\p_{u^j})
\\&
=-g(S^f(\p_{u^i},\p_{u^j}),(\p_{t}f)^\bot)
\\
g(\nabla^g_{\p_{u^i}}(\p_{t} f)^\top,Tf.\p_{u^j})
&=  (f^*g)(\nabla^{f^*g}_{\p_{u^i}}(\p_{t}f)^\top,\p_{u^j}) +0
\\
&= (f^*g)(\nabla^{f^*g}_{(\p_{t}f)^\top}\p_{u^i} + \on{Tor}
-[(\p_{t}f)^\top,\p_{u^i}],\p_{u^j}),
\\
\p_t(f^*g)_{ij}&= -2g(S^f(\p_{u^i},\p_{u^j}),(\p_{t}f)^\bot)
+(\mathcal L_{(\p_{t}f)^\top}(f^*g))(\p_{u^i},\p_{u^j})
\end{align*}
This proves the first formula. 
For the second one note that 
\begin{align*}
\tfrac12\on{Tr}((f^*g)\i\mathcal L_{h^\top}(f^*g))\on{vol}(f^*g)
= \mathcal L_{h^\top}(\on{vol}(f^*g)) =
\on{div}^{f^*g}(h^\top)\on{vol}(f^*g).
\qed\end{align*}
\end{demo}

\section{Metrics on spaces of mappings} \label{nmb:3}

\subsection{The metric $G^A$}\label{nmb:3.1}
Let $h,k\in C^\infty_f(M,TN)$ be two tangent vectors with foot point 
$f\in \on{Imm}(M,N)$, i.e., vector fields along $f$. Let the induced 
volume density be $\on{vol}(f^*g)$. We consider the following weak 
Riemannian metric on $\on{Imm}(M,N)$, for a constant $A\ge 0$:
\begin{equation*}
G^A_f(h,k) :=
\int_{M}\left(1+A\|\on{Tr}^{f^*g}(S^f)\|_{g^{N(f)}}^2\right)g(h,k)\on{vol}(f^*g)
\end{equation*}
where $\on{Tr}^{f^*g}(S^f) \in N(f)$ is the mean curvature, a section of the
normal bundle, and $\|\on{Tr}^{f^*g}(S^f)\|_{g^{N(f)}}$ is its norm.
The metric $G^A$ is invariant for the action of $\on{Diff}(M)$. This makes 
the map $\pi:\on{Imm}(M,N) \rightarrow B_i(M,N)$ into a {\it Riemannian 
submersion} (off the singularities of $B_i(M,N)$).

Now we can determine the bundle $\mathcal N\to \on{Imm}(M,N)$ 
of tangent vectors which are normal to the $\on{Diff}(M)$-orbits.
The tangent vectors to the orbits are $T_f(f\o\on{Diff}(M)) = 
\{Tf.\xi:\xi\in\X(M)\}$. Inserting this for $k$
into the expression of the metric $G$ we see that 
\begin{align*}
\mathcal N_f & =\{h\in C^\infty(M,TN): g(h,Tf)=0\}
\\
&=\Ga(N(f)),
\end{align*}
the space of sections of the normal bundle. This is independent of $A$.

A tangent vector $h\in
T_f\on{Imm}(M,N)=C^\infty_f(M,TN)=\Ga(f^*TN)$
has an orthonormal decomposition 
\begin{align*}
h&=h^\top+h^\bot \in T_f(f\o\on{Diff}^+(M)) \oplus \mathcal N_f
\notag\end{align*}
into smooth tangential and normal components.

Since the Riemannian metric $G^A$ on $\on{Imm}(M, N)$ is invariant
under the action of $\on{Diff}(M)$ it induces a metric on the quotient
$B_i(M, N)$ as follows. For any $F_0, F_1 \in B_i$, consider all
liftings $f_0, f_1 \in \on{Imm}$ such that $\pi(f_0)=F_0, \pi(f_1)=F_1$ and
all smooth curves $t\mapsto f(t,\quad)$ in $\on{Imm}(M,N)$
with $f(0,\cdot)=f_0$ and $f(1,\cdot)=f_1$. 
Since the metric $G^A$ is invariant under the action of $\on{Diff}(M)$
the arc-length of the curve 
$t\mapsto \pi(f(t,\cdot))$ in $B_i(M,N)$ is given by 
\begin{align*}
L^{\text{hor}}_{G^A}(f) &:= L_{G^A}(\pi(f(t,\cdot)))
=\\&
= \int_0^1 \sqrt{G^A_{\pi(f)}(T_f\pi.f_t,T_f\pi.f_t)}\,dt
 = \int_0^1 \sqrt{G^A_{f}(f_t^\bot,f_t^\bot)}\,dt
=\\&
=\int_0^1\Bigl(\int_{M}(1+A\|\on{Tr}^{f^*g}(S^f)\|_{g^{N(f)}}^2)
g(f_t^\bot,f_t^\bot)\on{vol}(f^*g)
\Bigr)^{\tfrac12}dt
\end{align*}

In fact the last computation only makes sense on $B_{i,f}(M,N)$ but 
we take it as a motivation.
The metric on $B_i(M, N)$ is defined by taking the infimum of this 
over all paths $f$ (and all lifts $f_0, f_1)$:
\begin{displaymath}
\on{dist}^{B_i}_{G^A}(F_1,F_2) = 
\inf_f L^{\text{hor}}_{G^A}(f).
\end{displaymath}

\begin{thm}\label{nmb:3.2}
Let $A=0$. For $f_0,f_1\in\on{Imm}(M,N)$ there exists always
a path $t\mapsto f(t,\cdot)$ in $\on{Imm}(M,N)$ with 
$f(0,\cdot)=f_0$ and $\pi(f(1,\cdot))=\pi(f_1)$ 
such that $L^{\text{hor}}_{G^0}(f)$ is arbitrarily small.
\end{thm}

\begin{demo}{Proof}
Take a path $f(t,\th)$ in $\on{Imm}(M,N)$ from $f_0$ to $f_1$ and
make it horizontal using \ref{nmb:2.4} so that
that $g(f_t,\D f)=0$; this forces a reparametrization on
$f_1$. 

Let $\al:M\to [0,1]$ be a surjective Morse function whose
singular values are all contained in the set
$\{\frac{k}{2N}:0\le k\le 2N\}$ for some integer $N$. We shall
use integers $n$ below and we shall use only multiples of $N$. 

Then the level sets
$M_r:=\{x\in M:\al(x)=r\}$ are of Lebesque measure 0.
We shall also need the slices 
$M_{r_1,r_2}:=\{x\in M:r_1\le\al(x)\le r_2\}$.
Since $M$ is compact there exists a constant $C$ such that the
following estimate holds uniformly in $t$:
$$
\int_{M_{r_1,r_2}} \on{vol}(f(t,\quad)^*g) 
\le C(r_2-r_1)\int_{M} \on{vol}(f(t,\quad)^*g) 
$$

Let $\tilde f(t,x)=f(\ph(t,\al(x)),x)$ where 
$\ph:[0,1]\x[0,1]\to[0,1]$ is given as in \cite{M98},~3.10 (which
also contains a figure illustrating the construction) by
\begin{displaymath}
\ph(t,\al)=
\begin{cases}
  2t(2n\al-2k) &\text{ for }\; 0 \le t\le 1/2,\;
    \tfrac{2k}{2n}\le \al \le\tfrac{2k+1}{2n}
  \\
  2t(2k+2-2n\al) &\text{ for }\; 0 \le t\le 1/2,\;
    \tfrac{2k+1}{2n}\le \al \le\tfrac{2k+2}{2n}
  \\
  2t-1+2(1-t)(2n\al-2k) &\text{ for }\; 1/2 \le t\le 1,\;
    \tfrac{2k}{2n}\le \al \le\tfrac{2k+1}{2n}
  \\
  2t-1+2(1-t)(2k+2-2n\al) &\text{ for }\; 1/2 \le t\le 1,\;
    \tfrac{2k+1}{2n}\le \al \le\tfrac{2k+2}{2n}.
\end{cases}
\end{displaymath}
Then we get $\D \tilde f = \ph_\al.d\al.f_t + \D f$ and 
$\tilde f_t = \ph_t. f_t$ where
\begin{displaymath}
\ph_\al=\begin{cases} +4nt \\ -4nt \\+4n(1-t) \\ -4n(1-t) \end{cases},\qquad 
\ph_t=\begin{cases} 4n\al-4k \\ 4k+4-4n\al \\2-4n\al+4k\\-(2-4n\al+4k)\end{cases}.
\end{displaymath}
We use horizontality $g(f_t,\D f)=0$ 
to determine $\tilde f_t^\bot=\tilde f_t + \D\tilde f(X)$ where
$X\in TM$ satisfies
$0=g(\tilde f_t + \D\tilde f(X),\D\tilde f(\xi))$ for all
$\xi\in TM$. We also use
$$
d\al(\xi)=f^*g(\on{grad}^{f^*g}\al,\xi)
  =g(\D f(\on{grad}^{f^*g}\al),\D f(\xi))
$$
and get
\begin{align*}
0 &= g(\tilde f_t + \D\tilde f(X),\D\tilde f(\xi))
\\&
= g\Bigl(\ph_tf_t + \ph_\al d\al(X) f_t + \D f(X), 
                \ph_\al d\al(\xi) f_t +\D f(\xi)\Bigr)
\\&
=\ph_t.\ph_\al.(f^*g)(\on{grad}^{f^*g}\al,\xi)\|f_t\|_g^2 
+\\&\quad
+\ph_\al^2.(f^*g)(\on{grad}^{f^*g}\al,X).(f^*g)(\on{grad}^{f^*g}\al,\xi)\|f_t\|_g^2
+g(\D f(X),\D f(\xi))
\\&
=(\ph_t.\ph_\al+\ph_\al^2.(f^*g)(\on{grad}^{f^*g}\al,X))\|f_t\|_g^2
(f^*g)(\on{grad}^{f^*g}\al,\xi) + (f^*g)(X,\xi)
\end{align*}
This implies that $X=\la\on{grad}^{f^*g}\al$ for a function $\la$
and in fact we get  
\begin{align*}
\tilde f_t^\bot 
= \frac{\ph_t}
   {1+\ph_\al^2\|d\al\|_{f^*g}^2\|f_t\|_g^2}
   f_t
-\frac{\ph_t \ph_\al\|f_t\|_g^2}
   {1+\ph_\al^2\|d\al\|_{f^*g}^2\|f_t\|_g^2}
  \D f(\on{grad}^{f^*g}\al)
\end{align*}
and 
$$
\|\tilde f_t\|_g^2 
= \frac{\ph_t^2\|f_t\|_g^2}
      {1+\ph_\al^2\|d\al\|_{f^*g}^2\|f_t\|_g^2}
$$
From $\D \tilde f = \ph_\al.d\al.f_t + \D f$ and 
$g(f_t,\D f)=0$ we get for the volume form 
$$
\on{vol}(\tilde f^*g) =
\sqrt{1+\ph_\al^2\,\|d\al\|_{f^*g}^2\|f_t\|_g^2}\,\on{vol}(f^*g).
$$
For the horizontal length we get
\begin{align*}
&L^{\text{hor}}(\tilde f)  
=\int_0^1\Bigl(\int_M
  \|\tilde f_t^\bot\|_g^2 \on{vol}(\tilde f^*g)\Bigr)^{\tfrac12}dt
=\\&
=\int_0^1\Bigl(\int_M
  \frac{\ph_t^2\|f_t\|_g^2}
      {\sqrt{1+\ph_\al^2\|d\al\|_{f^*g}^2\|f_t\|_g^2}}
      \on{vol}(f^*g)\Bigr)^{\tfrac12}dt
=\\&
=\int_0^{\frac12}\Biggl(\sum_{k=0}^{n-1}\Bigl(
\int_{M_{\tfrac{2k}{2n},\tfrac{2k+1}{2n}}} 
  \frac{(4n\al-4k)^2\|f_t\|_g^2}
      {\sqrt{1+(4nt)^2\|d\al\|_{f^*g}^2\|f_t\|_g^2}}
  \on{vol}(f^*g)
+\\&\qquad\qquad
+\int_{M_{\tfrac{2k+1}{2n},\tfrac{2k+2}{2n}}} 
  \frac{(4k+4-4n\al)^2\|f_t\|_g^2}
      {\sqrt{1+(4nt)^2\|d\al\|_{f^*g}^2\|f_t\|_g^2}}
      \on{vol}(f^*g)
\Bigr)\Biggr)^{\tfrac12}dt
+\\&\quad
+\int_{\frac12}^1\Biggl(\sum_{k=0}^{n-1}\Bigl(
\int_{M_{\tfrac{2k}{2n},\tfrac{2k+1}{2n}}} 
  \frac{(2-4n\al+4k)^2\|f_t\|_g^2}
      {\sqrt{1+(4n(1-t))^2\|d\al\|_{f^*g}^2\|f_t\|_g^2}}
  \on{vol}(f^*g)
+\\&\qquad\qquad\quad
+\int_{M_{\tfrac{2k+1}{2n},\tfrac{2k+2}{2n}}} 
  \frac{(2-4n\al+4k)^2\|f_t\|_g^2}
      {\sqrt{1+(4n(1-t))^2\|d\al\|_{f^*g}^2\|f_t\|_g^2}}
  \on{vol}(f^*g)
\Bigr)\Biggr)^{\tfrac12}dt
\end{align*}
Let $\ep>0$.
The function 
$(t,x)\mapsto \|f_t(\ph(t,\al(x)),x)\|_g^2$ is 
uniformly bounded. 
On $M_{\tfrac{2k}{2n},\tfrac{2k+1}{2n}}$ the function
$4n\al-4k$ has values in $[0,2]$. Choose disjoint geodesic
balls centered at the finitely many singular values of the
Morse function $\al$ of total $f^*g$-volume $<\ep$. Restricted
to the union $M_{\text{sing}}$ of these balls the integral
above is $O(1)\ep$. So we have to estimate the integrals on the
complement $\tilde M=M\setminus M_{\text{sing}}$ where the
function $\|d\al\|_{f^*g}$ is uniformly bounded from below by
$\et>0$.

Let us estimate one of the sums above. We use the fact that
the singular points of the Morse function $\al$ lie all on the
boundaries of the sets $\tilde M_{\tfrac{2k}{2n},\tfrac{2k+1}{2n}}$
so that we can transform the integrals as follows:
\begin{align*}
&\sum_{k=0}^{n-1}
\int_{\tilde M_{\tfrac{2k}{2n},\tfrac{2k+1}{2n}}} 
  \frac{(4n\al-4k)^2\|f_t\|_g^2}
      {\sqrt{1+(4nt)^2\|d\al\|_{f^*g}^2\|f_t\|_g^2}}
  \on{vol}(f^*g)=
\\&=\sum_{k=0}^{n-1}
\int_{\tfrac{2k}{2n}}^{\tfrac{2k+1}{2n}}\int_{\tilde M_{r}} 
  \frac{(4nr-4k)^2\|f_t\|_g^2}
      {\sqrt{1+(4nt)^2\|d\al\|_{f^*g}^2\|f_t\|_g^2}}
  \frac{\on{vol}(i_r^*f^*g)}{\|d\al\|_{f^*g}}\;dr
\end{align*}
We estimate this sum of integrals:
Consider first the set of all $(t,r,x\in M_r)$ such that 
$|f_t(\ph(t,r),x)|<\ep$. There we estimate by 
$$
O(1).n.16n^2.\ep^2.(r^3/3)|_{r=0}^{r=1/2n}=O(\ep).
$$
On the complementary set where $|f_t(\ph(t,r),x)|\ge\ep$ we
estimate by
$$
O(1).n.16n^2.\frac1{4nt\et^2\ep}(r^3/3)|_{r=0}^{r=1/2n}
  =O(\frac1{nt\et^2\ep})
$$
which goes to 0 if $n$ is large enough.
The other sums of integrals can be estimated similarly, thus 
$L^{\text{hor}}(\tilde f)$ goes to 0 for $n\to \infty$. It is clear that
one can approximate $\ph$ by a smooth function whithout changing the
estimates essentially. 
\qed\end{demo}

\subsection{A Lipschitz bound for the volume in $G^A$}\label{nmb:3.3}
We apply the Cauchy-Schwarz inequality to the 
derivative \ref{nmb:2.7} of the volume $\on{Vol}^g(f)$ along a
curve $t\mapsto f(t,\quad)\in\on{Imm}(M,N)$:
\begin{align*}
\p_t\on{Vol}^g(f)&
=\p_t \int_M \on{vol}^g(f(t,\quad))
= \int_M d\on{vol}^g(f)(\p_t f)
\\&
= -\int_M \on{Tr}^{f^*g}(g(S^f,f_t^\bot))\on{vol}(f^*g)
\le   \Bigl|\int_M \on{Tr}^{f^*g}(g(S^f,f_t^\bot))\on{vol}(f^*g)\Bigr|
\\&
\le \Bigl(\int_M 1^2\on{vol}(f^*g)\Bigr)^{\frac12}
 \Bigl(\int_M \on{Tr}^{f^*g}(g(S^f,f_t^\bot))^2\on{vol}(f^*g)
 \Bigr)^{\frac12}
\\&
\le \on{Vol}^g(f)^{\frac12}\frac1{\sqrt{A}}
 \Bigl(\int_M  (1+A\|\on{Tr}^{f^*g}(S^f)\|_{g^{N(f)}}^2)\,
                  g(f_t^\bot,f_t^\bot)\on{vol}(f^*g)
 \Bigr) ^{\frac12}
\end{align*}
Thus 
\begin{multline*}
\p_t(\sqrt{\on{Vol}^g(f)})=\frac{\p_t\on{Vol}^g(f)}{2\sqrt{\on{Vol}^g(f)}}
\le\\
\le \frac1{2\sqrt{A}}
 \Bigl(\int_M  (1+A\|\on{Tr}^{f^*g}(S^f)\|_{g^{N(f)}}^2)\,
                  g(f_t^\bot,f_t^\bot)\on{vol}(f^*g)
 \Bigr) ^{\frac12}
\end{multline*}
and by using \thetag{\ref{nmb:3.1}} we get
\begin{align*}
\sqrt{\on{Vol}^g(f_1)}&-\sqrt{\on{Vol}^g(f_0)}
=\int_0^1\p_t(\sqrt{\on{Vol}^g(f)})\,dt
\\&
\le \frac1{2\sqrt{A}}
 \int_0^1\Bigl(\int_M  (1+A\|\on{Tr}^{f^*g}(S^f)\|_{g^{N(f)}}^2)\,
                  g(f_t^\bot,f_t^\bot)\on{vol}(f^*g)
 \Bigr) ^{\frac12}dt
\\&
= \frac1{2\sqrt{A}}L^{\text{hor}}_{G^A}(f).
\end{align*}
If we take the infimum over all curves connecting $f_0$ with the
$\on{Diff}(M)$-orbit through $f_1$ we get:

\begin{prop*}{\bf Lipschitz continuity of 
$\sqrt{\on{Vol}^g}:B_i(M,N)\to \mathbb R_{\ge0}$.}
For $F_0$ and $F_1$ in $B_i(M,N)=\on{Imm}(M,N)/\on{Diff}(M)$ we have for $A>0$: 
\begin{equation*}
\sqrt{\on{Vol}^g(F_1)}-\sqrt{\on{Vol}^g(F_0)}\le
\frac1{2\sqrt{A}}\on{dist}^{B_i}_{G^A}(F_1,F_2).
\end{equation*}
\end{prop*}

\subsection{Bounding the area swept by a path in $B_i$}\label{nmb:3.4}
We want to bound the area swept out by a curve starting from $F_0$ to any
immersed submanifold $F_1$ nearby in our metric. 
First we use the Cauchy-Schwarz inequality 
in the Hilbert space $L^2(M,\on{vol}(f(t,\quad)^*g))$ to get  
\begin{multline*}
\int_{M}1.\|f_t\|_g\on{vol}(f^*g) = \langle 1,\|f_t\|_g\rangle_{L^2}\le
\\
\le \|1\|_{L^2}\|c_t\|_{L^2}
= \Bigl(\int_{M}\on{vol}(f^*g)\Bigr)^{\tfrac12}
\Bigl(\int_{M}|f_t|\on{vol}(f^*g)\Bigr)^{\tfrac12}.
\end{multline*}
Now we assume that the variation $f(t,x)$ is horizontal, so that 
$g(f_t,\D f)=0$. Then
$L_{G^A}(f)=L^{\text{hor}}_{G^A}(f)$.
We use this inequality and then the intermediate value theorem of integral
calculus to obtain 
\begin{align*}
L^{\text{hor}}_{G^A}(f) &= L_{G^A}(f) 
= \int_0^1 \sqrt{G^A_{f}(f_t,f_t)}\,dt
\\&
=\int_0^1\Bigl(\int_{M}(1+A\|\on{Tr}^{f^*}(S^f)\|_{f^*g}^2)
  \|f_t\|^2\on{vol}(f^*g)\Bigr)^{\tfrac12}dt
\\&
\ge\int_0^1\Bigl(\int_{M}\|f_t\|^2\on{vol}(f^*g)\Bigr)^{\tfrac12}dt
\\&
\ge\int_0^1 \Bigl(\int_{M}\on{vol}(f(t,\quad)^*g)\Bigr)^{-\tfrac12}
\int_{M}\|f_t(t,\quad)\|_g\on{vol}(f(t,\quad)^*g)\,dt
\\&
= \Bigl(\int_{M}\on{vol}(f(t_0,\quad)^*g)\Bigr)^{-\tfrac12}
\int_0^1\int_{M}\|f_t(t,\quad)\|_g\on{vol}(f(t,\quad)^*g)\,dt
\\&\qquad\qquad\text{  for some intermediate value $0\le t_0\le 1$, }
\\&
\ge \frac{1}{\sqrt{\on{Vol}^g(f(t_0,\quad))}}
\int_{[0,1]\x M}\on{vol}^{m+1}(f^*g)
\end{align*}

\begin{prop*}{\bf Area swept out bound.} 
If $f$ is any path from $F_0$ to $F_1$, then
\begin{equation*}
\begin{pmatrix} (m+1)-\text{volume of the region swept} \\
\text{out by the variation $f$}
  \end{pmatrix} \le \max_{t} \sqrt{\on{Vol}^g(f(t,\quad))} \cdot
L^{\text{hor}}_{G^A}(f). 
\end{equation*}
\end{prop*}
Together with the Lipschitz continuity \ref{nmb:3.3} this shows
that the geodesic distance $\inf L^{B_i}_{G^A}$ separates points, at least
in the base space $B(M,N)$ of embeddings.

\subsection{Horizontal energy of a path as anisotropic
volume}\label{nmb:3.5}
We consider a path $t\mapsto f(t,\quad)$ in $\on{Imm}(M,N)$. It
projects to a path $\pi\o f$ in $B_i$ whose energy is:
\begin{multline*}
E_{G^A}(\pi\o f) 
= \tfrac12 \int_a^b G^A_{\pi(f)}(T\pi.f_t,T\pi.f_t)\,dt
= \tfrac12 \int_a^b G^A_{f}(f_t^\bot,f_t^\bot)\,dt
=\\
= \tfrac12 \int_a^b 
  \int_{M}(1+A\|\on{Tr}^{f^*g}(S^f)\|_{g^{N(f)}}^2)
  g(f_t^\bot,f_t^\bot)\on{vol}(f^*g)\,dt.
\end{multline*}
We now consider the graph 
$\ga_f:[a,b]\x M\ni (t,x)\mapsto (t,f(t,x))\in [a,b]\x N$ 
of the path $f$ and its image $\Ga_f$, an immersed submanifold with
boundary of $\mathbb R\x N$. We want to describe the horizontal
energy as a functional on the space of immersed submanifolds
with fixed boundary, remembering the fibration of 
$\on{pr}_1:\mathbb R\x N\to \mathbb R$. We get:
\begin{align*}
&E_{G^A}(\pi\o f) =
\\&
= \tfrac12 
  \int_{[a,b]\x M}
  \left(1+A\|\on{Tr}^{f^*g}(S^f)\|_{g^{N(f)}}^2\right)
  \frac{\|f_t^\bot\|^2}{\sqrt{1+\|f_t^\bot\|_g^2}}
  \on{vol}(\ga_f^*(dt^2+g))
\end{align*}
Now $\|f_t^\bot\|_g$ depends only on the graph $\Ga_f$ and on
the fibration over time, since any reparameterization of
$\Ga_f$ which respects the fibration over time is of the form 
$(t,x)\mapsto(t,f(t,\ph(t,x)))$ for some path $\ph$ in $\on{Diff}(M)$
starting at the identity, and
$(\p_t|_0 f(t,\ph(t,x)))^\bot=f_t^\bot$. So the above
expression is intrinsic for the graph $\Ga_f$ and the
fibration. In order to find a geodesic from the shape 
$\pi(f(a,\quad))$ to the shape $\pi(f(b,\quad))$ one has to
find an immersed surface which is a critical point for the
functional $E_{G^A}$ above. This is a Plateau-problem with
anisotropic volume. 

\section{The geodesic equation and the curvature on
$B_i$}\label{nmb:4}

\subsection{The geodesic equation of $G^0$ in $\on{Imm}(M,N)$
}\label{nmb:4.1}
The energy of a curve $t\mapsto f(t,\quad)$ in $\on{Imm}(M,N)$
for $G^0$ is 
$$
E_{G^0}(f) = \tfrac12\int_a^b \int_M g(f_t,f_t)\on{vol}(f^*g).
$$
{\bf The geodesic equation for $G^0$}
\begin{equation}\boxed{\quad
{\begin{aligned}
\nabla^g_{\p_t} f_t &+ \on{div}^{f^*g}(f_t^\top)f_t 
- g(f_t^\bot,\on{Tr}^{f^*g}(S^f))f_t
+\\&
+ \tfrac12 \D f.\on{grad}^{f^*g}(\|f_t\|_g^2) 
+ \tfrac12 \|f_t\|_g^2 \on{Tr}^{f^*g}(S^f) =0
\end{aligned}}\quad
}
\tag{1}
\end{equation}

\begin{demo}{Proof} A different proof is in \cite{Ka},~2.2. 
For a function $a$ on $M$ we shall use 
\begin{align*}
\int_M a\on{div}^{f^*g}(X)\on{vol}(f^*g) 
&= \int_M a\mathcal L_X(\on{vol}(f^*g)) 
\\&
= \int_M \mathcal L_X(a\on{vol}(f^*g)) -    
 \int_M \mathcal L_X(a)\on{vol}(f^*g) 
\\&
=- \int_M (f^*g)(\on{grad}^{f^*g}(a),X)\on{vol}(f^*g)
\end{align*}
in calculating the first variation of the energy with fixed
ends:
\begin{align*}
\p_s E_{G^0}(f) &=
\tfrac12\int_a^b\int_M \Bigl(\p_s g(f_t,f_t)\,\on{vol}(f^*g) 
  +g(f_t,f_t)\,\p_s\on{vol}(f^*g)\Bigr)dt
\\&
= \int_a^b\int_M \Bigl(g(\nabla^g_{\p_s}f_t,f_t)\,\on{vol}(f^*g) 
  +\tfrac12 \|f_t\|_g^2\on{div}^{f^*g}(f_s^\top)\on{vol}(f^*g)
\\&\qquad\qquad\qquad\qquad\qquad\qquad
  -\tfrac12 \|f_t\|_g^2 g(f_s^\bot,\on{Tr}^{f^*g}(S^f))\on{vol}(f^*g)
  \Bigr)dt
\end{align*}
For the first summand we have:
\begin{align*}
&\int_a^b\int_M g(\nabla^g_{\p_s}f_t,f_t)\,\on{vol}(f^*g)\,dt
=\int_a^b\int_M g(\nabla^g_{\p_t}f_s,f_t)\,\on{vol}(f^*g)\,dt
\\&
=\int_a^b\int_M (\p_t g(f_s,f_t)
  -g(f_s,\nabla^g_{\p_t}f_t))\on{vol}(f^*g)\,dt
\\&
=-\int_a^b\int_M g(f_s,f_t)\p_t\on{vol}(f^*g)\,dt
-\int_a^b\int_M g(f_s,\nabla^g_{\p_t}f_t)\on{vol}(f^*g)\,dt
\\&
=\int_a^b\int_M \Bigl(-g(f_s,f_t)\on{div}^{f^*g}(f_t^\top) 
+g(f_s,f_t)g(f_t^\bot, \on{Tr}^{f^*g}(S^f))
-\\&\qquad\qquad\qquad\qquad\qquad\qquad\qquad\qquad\qquad\qquad
-g(f_s,\nabla^g_{\p_t}f_t)\Bigr)\on{vol}(f^*g)\,dt
\end{align*}
The second summand yields:
\begin{align*}
\int_a^b\int_M
\tfrac12 &\|f_t\|_g^2 \on{div}^{f^*g}(f_s^\top)\on{vol}(f^*g)\,dt 
\\&
=-\int_a^b\int_M \tfrac12 (f^*g)(f_s^\top,\on{grad}^{f^*g}(\|f_t\|_g^2))\on{vol}(f^*g)\,dt
\\&
=-\int_a^b\int_M \tfrac12 g(f_s,\D f.\on{grad}^{f^*g}(\|f_t\|_g^2))\on{vol}(f^*g)\,dt
\end{align*}
Thus the first variation $\p_s E_{G^0}(f)$ is:
\begin{align*}
\int_a^b\int_M 
g\Bigl(f_s, &-\nabla^g_{\p_t}f_t
             -\on{div}^{f^*g}(f_t^\top)f_t 
             +g(f_t^\bot, \on{Tr}^{f^*g}(S^f))f_t
\\&
             -\tfrac12 \D f.\on{grad}^{f^*g}(\|f_t\|_g^2)
             -\tfrac12 \|f_t\|_g^2 g(f_s^\bot,\on{Tr}^{f^*g}(S^f))
  \Bigr)\on{vol}(f^*g)\,dt
\qed\end{align*}
\end{demo}

\subsection{Geodesics for $G^0$ in $B_i(M,N)$
}\label{nmb:4.2}
We restrict to geodesics $t\mapsto f(t,\quad)$ in $\on{Imm}(M,N)$ which are
horizontal: $g(f_t,\D f)=0$. Then $f_t^\top=0$ and
$f_t=f_t^\bot$, so equation \ref{nmb:4.1}.\thetag{1} becomes
$$
\nabla^g_{\p_t} f_t
- g(f_t,\on{Tr}^{f^*g}(S^f))f_t
+ \tfrac12 \D f.\on{grad}^{f^*g}(\|f_t\|_g^2) 
+ \tfrac12 \|f_t\|_g^2 \on{Tr}^{f^*g}(S^f) =0.
$$
It splits into a vertical (tangential) part
$$
-\D f.(\nabla_{\p_t}f_t)^\top +\tfrac12
\D f.\on{grad}^{f^*g}(\|f_t\|_g^2) =0
$$
which vanishes identically since  
\begin{multline*}
(f^*g)(\on{grad}^{f^*g}(\|f_t\|^2_g),X) = X(g(f_t,f_t)) =
2g(\nabla_Xf_t,f_t) = 2g(\nabla^g_{\p_t}\D f.X,f_t) 
\\
= 2\p_t g(\D f.X,f_t) - 2g(\D f.X,\nabla^g_{\p_t}f_t)
=  - 2g(\D f.X,\nabla^g_{\p_t}f_t),
\end{multline*}
and a horizontal (normal) part which is the geodesic equation in $B_i$:
\begin{equation}\boxed{
\nabla^{N(f)}_{\p_t} f_t
- g(f_t,\on{Tr}^{f^*g}(S^f))f_t
+ \tfrac12 \|f_t\|_g^2 \on{Tr}^{f^*g}(S^f) = 0,
\quad g(\D f,f_t)=0}
\tag{1}
\end{equation}

\subsection{The induced metric of $G^0$ in $B_i(M,N)$ in a
chart
}\label{nmb:4.3}
Let $f_0:M\to N$ be a fixed immersion which will be the `center' of our chart. 
Let $N(f_0)\subset f_0^*TN$ be the normal bundle to $f_0$. Let 
$\exp^g: N(f_0) \rightarrow N$ be the exponential map for the metric $g$ 
and let $V\subset N(f_0)$ be a neighborhood of the 0 section 
on which the exponential map is an immersion. Consider the mapping 
\begin{align}
&\ps=\ps_{f_0}:\Ga(V)\to \on{Imm}(M,N),\quad, \ps(\Ga(V))=:\mathcal Q(f_0),
\tag{1}\\&
\ps(a)(x) = \exp^g(a(x))=\exp^g_{f_0(x)}(a(x)).
\notag\end{align}
The inverse (on its image) of $\pi\o\ps_f:\Ga(V)\to B_i(M,N)$
is a smooth chart on $B_i(M,N)$. Our goal is to calculate the induced metric
on this chart, that is 
$$((\pi \circ \ps_{f_0})^* G^0_a )(b_1,b_2)$$ 
for any $a \in \Ga(V), b_1,b_2 \in \Ga(N(f_0))$. This will enable us to 
calculate the sectional curvatures of $B_i$.

We shall fix the section $a$ and work with the ray of points $t.a$ in this
chart. Everything will revolve around the map:
$$ f(t,x) = \ps(t.a)(x) = \exp^g(t.a(x)).$$
We shall also use a fixed chart $(M\supset U \East{u}{} \mathbb
R^m)$ on $M$ with $\p_i=\p/\p u^i$.
Then $x\mapsto (t\mapsto f(t,x))=\exp^g_{f_0(x)}(t.a(x))$ 
is a variation consisting entirely of geodesics, thus:
\begin{align*}
t\mapsto \p_i f(t,x)&=Tf.\p_i=:Z_i(t,x,a)
 \text{ \it is the Jacobi field along } t\mapsto f(t,x)
 \text{ \it with }
\\
Z_i(0,x,a) &= \p_i|_x \exp_{f_0(x)}(0) = \p_i|_x f_0 =
T_xf_0.\p_i|_x,
\tag{2}\\
(\nabla^g_{\p_t}Z_i)(0,x) &= (\nabla^g_{\p_t}Tf.\p_i)(0,x) =
(\nabla^g_{\p_i}Tf.\p_t)(0,x) =
\\&
=\nabla^g_{\p_i}(\p_t|_0\exp_{f_0(x)}(t.a(x))) =
(\nabla^g_{\p_i}a)(x).
\end{align*}
Then the pullback metric is given by
\begin{align}
f^*g&=\ps(ta)^*g = g(\D f,\D f) 
= \sum_{i,j=1}^m g(Tf.\p_i,Tf,\p_j)\; du^i\otimes du^j
\notag\\& 
= \sum_{i,j=1}^m g(Z_i,Z_j)\; du^i\otimes du^j.
\tag{3}\end{align}
The induced volume density is:
\begin{align}
\on{vol}(f^*g)&=\sqrt{\det(g(Z_i,Z_j))}|du^1\wedge \dots\wedge du^m|
\tag{4}\end{align}
Moreover we have for $a\in \Ga(V)$ and $b\in\Ga(N(f_0))$ 
\begin{align}
(T_{ta}\ps.tb)(x) 
&= \p_s|_0 \exp^g_{f_0(x)}(ta(x)+stb(x)))
\notag\\& 
= Y(t,x,a,b)
\quad\text{ for the Jacobi field }Y\text{ along }t\mapsto
f(t,x)\text{ with}
\tag{5}\\
Y(0,x,a,b) &= 0_{f_0(x)},
\notag\\ 
(\nabla^g_{\p_t}Y&(\quad,x,a,b))(0) 
= \nabla^g_{\p_t}\p_s \exp^g_{f_0(x)}(ta(x)+stb(x))
\notag\\& 
= \nabla^g_{\p_s}\p_t|_0 \exp^g_{f_0(x)}(ta(x)+stb(x))
\notag\\& 
= \nabla^g_{\p_s}(a(x)+s.b(x))|_{s=0} = b(x).
\notag\end{align}
Now we want to split $T_a\ps.b$ into vertical (tangential) and
horizontal parts with respect to the immersion
$\ps(ta)=f(t,\quad)$.
The tangential part has locally the form
\begin{align*}
&\D f.(T_{ta}\ps.tb)^\top 
= \sum_{i=1}^m c^i\; \D f.\p_i
= \sum_{i=1}^m c^i\; Z_i\quad\text{  where for all }j
\\&
g(Y,Z_j) 
= g\left(\sum_{i=1}^m c^i\; Z_i,Z_j \right) 
= \sum_{i=1}^m c^i\; (f^*g)_{ij},\quad
\\&
c^i= \sum_{j=1}^m (f^*g)^{ij} g(Y,Z_j).
\end{align*}
Thus the horizontal part is 
\begin{align}
(T_{ta}\ps.tb)^\bot = Y^\bot = Y-\sum_{i=1}^m c^i\; Z_i
= Y-\sum_{i,j=1}^m (f^*g)^{ij}\; g(Y,Z_j)\; Z_i 
\tag{6}\end{align}
Thus the induced metric on $B_i(M,N)$ has the following
expression in the chart $(\pi\o\ps_{f_0})\i$, where
$a\in\Ga(V)$ and $b_1,b_2\in\Ga(N(f_0))$:
\begin{align}
((\pi\o\ps_{f_0}&)^*G^0)_{ta}(b_1,b_2) 
= G^0_{\pi(\ps(ta))}(T_{ta}(\pi\o\ps)b_1,T_{ta}(\pi\o\ps)b_2)
\notag\\&
= G^0_{\ps(ta)}((T_{ta}\ps.b_1)^\bot,(T_{ta}\ps.b_2)^\bot)
\notag\\&
= \int_M g((T_{ta}\ps.b_1)^\bot,(T_{ta}\ps.b_2)^\bot)\;\on{vol}(f^*g)
\notag\\&
= \int_M  \frac1{t^2}g\Bigl(Y(b_1)-\sum_{i,j} (f^*g)^{ij}\, g(Y(b_1),Z_j)\,
  Z_i\;,\; Y(b_2)\Bigr) \;
\tag{7}\\&\qquad\qquad\qquad\qquad\qquad\qquad\qquad
  \sqrt{\det(g(Z_i,Z_j))}|du^1\wedge \dots\wedge du^m|
\notag\end{align}

\subsection{Expansion to order 2 of the induced metric of $G^0$ in $B_i(M,N)$ 
in a chart
}\label{nmb:4.4}
We use the setting of \ref{nmb:4.3}, the Einstein summation
convention, and the abbreviations $f_i:=\p_if_0=\p_{u^i}f_0$
and $\nabla^g_i:=\nabla^g_{\p_i}=\nabla^g_{\p_{u^i}}$.
We compute the expansion in $t$ up to order 2 of the metric
\ref{nmb:4.3}.\thetag{7}. 
Our method is to use the Jacobi equation 
$$
\nabla^g_{\p_t}\nabla^g_{\p_t}Y = R^g(\dot c,Y)\dot c
$$
which holds for any Jacobi field $Y$ along a geodesic $c$.
By \ref{nmb:2.6}.\thetag{2} we have:
\begin{equation}
g(\nabla^g_ia,f_j)=-(f_0^*g)(L^{f_0}_a(\p_i),\p_j)
  =-g(a,S^{f_0}(f_i,f_j))=g(a,S^{f_0}_{ij})
\tag{1}\end{equation}
We start by expanding the pullback metric
\ref{nmb:4.3}.\thetag{3} and its inverse:
\begin{align}
\p_t(f^*g)_{ij}&=
\p_t g(Z_i,Z_j) = g(\nabla^g_{\p_t}Z_i,Z_j) + 
g(Z_i,\nabla^g_{\p_t}Z_j) 
\notag\\
\p_t^2 g(Z_i,Z_j) &= g(\nabla^g_{\p_t}\nabla^g_{\p_t}Z_i,Z_j) + 
  2g(\nabla^g_{\p_t}Z_i,\nabla^g_{\p_t}Z_j) + 
g(Z_i,\nabla^g_{\p_t}\nabla^g_{\p_t}Z_j) 
\notag\\
(f^*g)_{ij}&=(f_0^*g)_{ij} 
  + t\bigl(g(\nabla^g_ia,f_j)+g(f_i,\nabla^g_ja)\bigr)
+\notag\\&\quad
  + \tfrac12{t^2}\bigl(g(R^g(a,f_i)a,f_j) + 
  2g(\nabla^g_ia,\nabla^g_ja) + 
  g(f_i,R^g(a,f_j)a)\bigr)
\notag\\&\quad
  +O(t^3)
\notag\\&
=(f_0^*g)_{ij} 
  -2 t (f_0^*g)(L^{f_0}_a(\p_i),\p_j)
\notag\\&\quad
  + t^2\bigl(g(R^g(a,f_i)a,f_j) + 
  g(\nabla^g_ia,\nabla^g_ja)\bigr) 
  +O(t^3)
\tag{2}
\end{align}
We expand now the volume form
$ \on{vol}(f^*g)=\sqrt{\det(g(Z_i,Z_j))}|du^1\wedge \dots\wedge du^m|$. 
The time derivative at 0 of the inverse of the pullback metric is:
\begin{align*}
\p_t(f^*g)^{ij}|_0 &= -(f_0^*g)^{ik}(\p_t|_0(f^*g)_{kl})(f_0^*g)^{lj}
= -(f_0^*g)^{ik}(f_0^*g)(a,S^{f_0}_{kl})(f_0^*g)^{lj}
\end{align*}
Therefore,
\begin{align*}
\p_t \sqrt{\det(g(Z_i,Z_j))} 
  &= \tfrac12 (f^*g)^{ij}\p_t(g(Z_i,Z_j))\sqrt{\det(g(Z_i,Z_j))}
\\
\p_t^2 \sqrt{\det(g(Z_i,Z_j))} 
&= \tfrac12 \p_t(f^*g)^{ij}\p_t(g(Z_i,Z_j))\sqrt{\det(g(Z_i,Z_j))}
\\&\quad 
+\tfrac12 (f^*g)^{ij}\p_t^2(g(Z_i,Z_j))\sqrt{\det(g(Z_i,Z_j))}
\\&\quad 
+\tfrac12 (f^*g)^{ij}\p_t(g(Z_i,Z_j))\p_t\sqrt{\det(g(Z_i,Z_j))},
\end{align*}
and 
\begin{align*}
\on{vol}(f^*g)&=\sqrt{\det(g(Z_i,Z_j))}|du^1\wedge \dots\wedge du^m|
\\&
= \Bigl(1 
- t\on{Tr}(L^{f_0}_a)
+ t^2\Bigl(-\on{Tr}(L^{f_0}_a\o L^{f_0}_a)
+\tfrac12(\on{Tr}(L^{f_0}_a))^2
\\&\qquad\qquad
+ \tfrac12(f^*g)^{ij}\bigl(g(R^g(a,f_i)a,f_j) + g(\nabla^g_ia,\nabla^g_ja)\bigr)
\Bigr)+O(t^3)
\Bigr)\on{vol}(f_0^*g)
\tag{3}\end{align*}
Moreover, by \ref{nmb:2.6}.\thetag{2} we may split $\nabla^g_i a =
-Tf_0.L^{f_0}_a(\p_i) + \nabla^{N(f_0)}_ia$ and we write $\nabla^\bot_ia$
for $\nabla^{N(f_0)}_ia$ shortly. Thus:
\begin{align*}
(f^*_0 g)^{ij} g(\nabla^g_i a, \nabla^g_j a) 
&= (f^*_0 g)^{ij}\left( g(Tf_0 . L^{f_0}_a(\p_i), Tf_0.L^{f_0}_a(\p_j)) + 
g(\nabla^\bot_i a, \nabla^\bot_j a)\right)\\
&= \on{Tr}(L^{f_0}_a \circ L^{f_0}_a) + 
(f^*_0 g)^{ij} g(\nabla^\bot_i a, \nabla^\bot_j a)
\end{align*}
and so that the tangential term above combines with the first $t^2$ term in 
the expansion of the volume, changing its coefficient from $-1$ to $-\tfrac12$.

Let us now expand
\begin{align*}
g((T_{ta}\ps. tb_1)^\bot,(T_{ta}\ps. tb_2)^\bot)
&= g\bigl(Y(b_1)- (f^*g)^{ij} g(Y(b_1),Z_j)\, Z_i\;,\; Y(b_2)\bigr)
\\&
= g(Y(b_1),Y(b_2))-(f^*g)^{ij}g(Y(b_1),Z_j)g(Z_i,Y(b_2)).
\end{align*}
We have: 
\begin{align*}
\p_t g(Y(b_1),Y(b_2))
&=g(\nabla^g_{\p_t}Y(b_1),Y(b_2))+g(Y(b_1),\nabla^g_{\p_t}Y(b_2))
\\
\p_t^2 g(Y(b_1),Y(b_2))
&=g(\nabla^g_{\p_t}\nabla^g_{\p_t}Y(b_1),Y(b_2))
+2g(\nabla^g_{\p_t}Y(b_1),\nabla^g_{\p_t}Y(b_2))
\\&\quad
+g(Y(b_1),\nabla^g_{\p_t}\nabla^g_{\p_t}Y(b_2))
\\&
=2g(R^g(a,Y(b_1))a,Y(b_2))
+2g(\nabla^g_{\p_t}Y(b_1),\nabla^g_{\p_t}Y(b_2))
\\
\p_t g(Y(b_1),Z_j)
&=g(\nabla^g_{\p_t}Y(b_1),Z_j)+g(Y(b_1),\nabla^g_{\p_t}Z_j)
\\
\p_t^2 g(Y(b_1),Z_j)
&=2g(R^g(a,Y(b_1))a,Z_j)
+2g(\nabla^g_{\p_t}Y(b_1),\nabla^g_{\p_t}Z_j)
\end{align*}
Note that: 
\begin{align*}
&Y(0,h)=0,\quad (\nabla^g_{\p_t}Y(h))(0)=h,\quad 
(\nabla^g_{\p_t}\nabla^g_{\p_t}Y(h))(0)=R^g(a,Y(0,h))a =0,
\\&
(\nabla^g_{\p_t}\nabla^g_{\p_t}\nabla^g_{\p_t}Y(h))(0)
= R^g(a,\nabla^g_{\p_t}Y(h)(0))a = R^g(a,h)a.
\end{align*}
Thus: 
\begin{align*}
&g(Y(b_1),Y(b_2)) 
-(f^*g)^{ij}g(Y(b_1),Z_j)g(Z_i,Y(b_2)) 
\\
&=  t^2 g(b_1,b_2) 
+t^4\bigl(\tfrac1{3} g(R^g(a,b_1)a,b_2)
-(f_0^*g)^{ij}g(b_1,\nabla^\bot_ja)g(\nabla^\bot_ia,b_2)\bigr) +O(t^5).
\end{align*}
The expansion of $G^0$ up to order 2 is thus:
\begin{align*}
((\pi\o\ps_{f_0}&)^*G^0)_{ta}(b_1,b_2) =
\notag\\&
= \int_M 
  \frac1{t^2}g\Bigl(Y(b_1)-\sum_{i,j} (f^*g)^{ij}\, g(Y(b_1),Z_j)\,
  Z_i\;,\; Y(b_2)\Bigr)\;\on{vol}(f^*g)
\\&
= \int_M\Bigl(g(b_1,b_2)\on{vol}(f_0^*g) 
-t\int_M g(b_1,b_2)\on{Tr}(L^{f_0}_a)\on{vol}(f_0^*g)
\\&\quad
+t^2 \int_M
\biggl( g(b_1,b_2) 
\Bigl(-\tfrac12 \on{Tr}(L^{f_0}_a\o L^{f_0}_a)
+\tfrac12 \on{Tr}(L^{f_0}_a)^2
\\&\qquad\qquad\qquad\qquad\quad
+ \tfrac12(f^*g)^{ij}g(R^g(a,f_i)a,f_j) 
+ \tfrac12(f^*g)^{ij}g(\nabla^\bot_ia,\nabla^\bot_ja) \Bigr)
\\&\qquad\qquad
+ \tfrac1{3} g(R^g(a,b_1)a,b_2)
-(f_0^*g)^{ij}g(b_1,\nabla^\bot_ja)g(\nabla^\bot_ia,b_2)\biggr)
\on{vol}(f_0^*g)
\\&+O(t^3)
\tag{4}\end{align*}

\subsection{ Computation of the sectional curvature in 
$B_i(M,N)$ at $f_0$}\label{nmb:4.5}
We use the following formula which is valid in a chart:
\begin{align*}
&2R_a(m,h,m,h)=2G^0_a(R_a(m,h)m,h) =
\\&
= -2d^2G^0(a)(m,h)(h,m)  +d^2G^0(a)(m,m)(h,h)  +d^2G^0(a)(h,h)(m,m)  
\\&\quad
-2G^0(\Ga(h,m),\Ga(m,h)) +2G^0(\Ga(m,m),\Ga(h,h)) 
\end{align*}
The sectional curvature at the two-dimensional subspace $P_a(m,h)$ of the
tangent space which is spanned by $m$ and $h$ is then given by:
\begin{equation*}
k_a(P(m,h)) = - \frac{G^0_a(R(m,h)m,h)}{\|m\|^2\|h\|^2-G^0_a(m,h)^2}.
\end{equation*}
We compute this directly for $a=0$.
From the expansion up to order 2 of $G^0_{ta}(b_1,b_2)$ in
\ref{nmb:4.4}.\thetag{4} we get
$$
dG^0(0)(a)(b_1,b_2) 
=-\int_M g(b_1,b_2) g(a,\on{Tr}^{f_0^*g}S^{f_0})\on{vol}(f_0^*g)
$$
and compute the Christoffel symbol:
\begin{align*}
-2G^0_0&(\Ga_0(a,b),c) = -dG^0(0)(c)(a,b) +dG^0(0)(a)(b,c)
+dG^0(0)(b)(c,a)
\\&
=\int_M\Bigl( g(a,b) g(c,\on{Tr}^{f_0^*g}(S^{f_0}))
- g(b,c) g(a,\on{Tr}^{f_0^*g}(S^{f_0}))
\\&\qquad\qquad
- g(c,a) g(b,\on{Tr}^{f_0^*g}(S^{f_0}))\Bigr)\on{vol}(f_0^*g)
\\&
=\int_M g\Bigl(c,\; g(a,b) \on{Tr}^{f_0^*g}(S^{f_0})
- \on{Tr}(L^{f_0}_a)b
- \on{Tr}(L^{f_0}_b)a\Bigr)\on{vol}(f_0^*g)
\\
\Ga_0&(a,b) = -\tfrac12 g(a,b) \on{Tr}^{f_0^*g}(S^{f_0})
+\tfrac12 \on{Tr}(L^{f_0}_a)b
+\tfrac12 \on{Tr}(L^{f_0}_b)a
\end{align*}
The expansion \ref{nmb:4.4}.\thetag{4} also gives:
\begin{align*}
&\tfrac1{2!}d^2G^0_{0}(a_1,a_2)(b_1,b_2) =
\\&
= \int_M
\biggl( g(b_1,b_2) 
\Bigl(-\tfrac12\on{Tr}(L^{f_0}_{a_1}\o L^{f_0}_{a_2})
+\tfrac12\on{Tr}(L^{f_0}_{a_1})\on{Tr}(L^{f_0}_{a_2})
\\&\qquad\qquad\qquad\qquad\qquad\quad
+ \tfrac12(f^*g)^{ij}g(R^g(a_1,f_i)a_2,f_j) +
\tfrac12(f^*g)^{ij}g(\nabla^\bot_ia_1,\nabla^\bot_ja_2)
\Bigr)
\\&\qquad\qquad\qquad
+ \tfrac1{6} g(R^g(a_1,b_1)a_2,b_2)
+ \tfrac1{6} g(R^g(a_2,b_1)a_1,b_2)
\\&\qquad\qquad\qquad
-\tfrac12(f_0^*g)^{ij}g(b_1,\nabla^\bot_ja_1)g(\nabla^\bot_ia_2,b_2)
\\&\qquad\qquad\qquad
-\tfrac12(f_0^*g)^{ij}g(b_1,\nabla^\bot_ja_2)g(\nabla^\bot_ia_1,b_2)
\biggr)\on{vol}(f_0^*g)
+O(t^3)
\end{align*}
Thus we have:
\begin{align*}
&-d^2G^0(0)(x,y)(y,x)+\tfrac12 d^2G^0(0)(x,x)(y,y)+\tfrac12 d^2G^0(0)(y,y)(x,x) =
\\&
= \int_M
\biggl(- 2g(y,x) 
\Bigl(-\tfrac12\on{Tr}(L^{f_0}_{x}\o L^{f_0}_{y})
+\tfrac12\on{Tr}(L^{f_0}_{x})\on{Tr}(L^{f_0}_{y})
\\&\qquad\qquad\qquad\qquad\qquad\quad
+ \tfrac12(f_0^*g)^{ij}g(R^g(x,f_i)y,f_j) 
+ \tfrac12(f_0^*g)^{ij}g(\nabla^\bot_ix,\nabla^\bot_jy)
\Bigr)
\\&\qquad\qquad\qquad 
+g(y,y) 
\Bigl(-\tfrac12\on{Tr}(L^{f_0}_{x}\o L^{f_0}_{x})
+\tfrac12\on{Tr}(L^{f_0}_{x})\on{Tr}(L^{f_0}_{x})
\\&\qquad\qquad\qquad\qquad\qquad\quad
+ \tfrac12(f_0^*g)^{ij}g(R^g(x,f_i)x,f_j) +
  \tfrac12(f_0^*g)^{ij}g(\nabla^\bot_ix,\nabla^\bot_jx)
\Bigr)
\\&\qquad\qquad\qquad   
+ g(x,x) 
\Bigl(-\tfrac12\on{Tr}(L^{f_0}_{y}\o L^{f_0}_{y})
+\tfrac12\on{Tr}(L^{f_0}_{y})\on{Tr}(L^{f_0}_{y})
\\&\qquad\qquad\qquad\qquad\qquad\quad
+ \tfrac12(f_0^*g)^{ij}g(R^g(y,f_i)y,f_j) +
  \tfrac12(f_0^*g)^{ij}g(\nabla^\bot_iy,\nabla^\bot_jy)
\Bigr)
\\&\qquad\qquad\qquad 
+ g(R^g(y,x)y,x) 
\\&\qquad\qquad\qquad 
+(f_0^*g)^{ij}\left(g(y,\nabla^\bot_jy)g(\nabla^\bot_ix,x)+
g(y,\nabla^\bot_jx)g(\nabla^\bot_iy,x)\right)
\\&\qquad\qquad\qquad 
-(f_0^*g)^{ij}\left(g(x,\nabla^\bot_jy)g(\nabla^\bot_iy,x)+
g(y,\nabla^\bot_jx)g(\nabla^\bot_ix,y)\right)
\biggr)\on{vol}(f_0^*g)
\end{align*}
For the second part of the curvature we have
\begin{align*}
&-G_0(\Ga_0(x,y),\Ga_0(y,x))+G_0(\Ga_0(y,y),\Ga_0(x,x)) =
\\&
=\tfrac14\!\int_M\!\Bigl(
(\|x\|_g^2\|y\|_g^2-g(x,y)^2)\|\on{Tr}^{f_0^*g}(S^f)\|_g^2  
-3\|\on{Tr}(L^{f_0}_x)y - \on{Tr}(L^{f_0}_y)x\|_g^2
\Bigr)\on{vol}(f_0^*g)
\end{align*}

To organize all these terms in the curvature tensor, note that they belong 
to {\it three} types: terms which involve the second fundamental 
form $L^{f_0}$, terms which involve the curvature tensor $R^g$ of $N$ and terms 
which involve the normal component of the covariant derivative $\nabla^\bot a$.
There are 3 of the first type, two of the second and the ones of the third 
can be organized neatly into two also. The final curvature tensor is the 
integral over $M$ of their sum. Here are the terms in detail:


\noindent
\thetag{1} {\it Terms involving the trace of products of $L$'s}. These are:
$$ -\tfrac12 \Big( g(y,y) \on{Tr}(L^{f_0}_x \circ L^{f_0}_x)
-2g(x,y)\on{Tr}(L^{f_0}_x \circ L^{f_0}_y) + 
g(x,x)\on{Tr}(L^{f_0}_y \circ L^{f_0}_y) \Big).$$
Note that $x$ and $y$ are sections of the normal bundle $N(f_0)$, so we may
define $x \wedge y$ to be the induced section of $\bigwedge^2 N(f_0)$. Then
the expression inside the parentheses is a positive semi-definite quadratic 
function of $x \wedge y$.  To see this, note a simple linear algebra fact --
that if $Q(a,b)$ is any positive semi-definite inner product on 
${\mathbb R}^n$, then 
\begin{align*}
\widetilde{Q}(a \wedge b, c \wedge d) &= \\
<a,c>&Q(b,d)-<a,d>Q(b,c)+<b,d>Q(a,c)-<b,c>Q(a,d) \\
\widetilde{Q}(a \wedge b, a \wedge b) &= 
\|a\|^2Q(b,b)-2<a,b>Q(a,b)+\|b\|^2Q(a,a)
\end{align*} 
is a positive semi-definite inner product on $\bigwedge^2 V$. In particular,
$\on{Tr}(L^{f_0}_x \circ L^{f_0}_y)$ is a positive semi-definite inner
product on the normal bundle, hence it defines a positive semi-definite
inner product $\widetilde{\on{Tr}}(L^{f_0} \circ L^{f_0})$ on 
$\bigwedge^2 N(f_0)$. Thus:
$$ \on{term (1)} 
  = -\tfrac12 \widetilde{\on{Tr}}(L^{f_0} \circ L^{f_0})(x \wedge y) \le 0.$$

\noindent
\thetag{2} {\it Terms involving trace of one $L$}.  We have terms both from 
the second and first derivatives of $G$, namely:
$$ \tfrac12 \Big( g(y,y) \on{Tr}(L^{f_0}_x)^2 
-2g(x,y)\on{Tr}(L^{f_0}_x) \on{Tr}(L^{f_0}_y) + 
g(x,x)\on{Tr}(L^{f_0}_y)^2 \Big).$$
and
$$ -\tfrac34 \| \on{Tr}(L^{f_0}_x)y-\on{Tr}(L^{f_0}_y)x\|^2_g$$
which are the same up to their coefficients. Their sum is:
$$\on{term (2)} = -\tfrac14 \| \on{Tr}(L^{f_0}_x)y-
\on{Tr}(L^{f_0}_y)x\|^2_g \le 0.$$
Note that this is a function of $x \wedge y$ also.

\noindent
\thetag{3} {\it The term involving the norm of the second fundamental form}.
Since $\|x\|_g^2 \|y\|^2_g - g(x,y)^2 = \|x \wedge y\|^2_g$, this term is just:
$$\on{term (3)} = +\tfrac14 \|x \wedge y\|_g^2 \|\on{Tr}^g(S^{f_0})\|_g^2 \ge 0.$$

\noindent
\thetag{4} {\it The curvature of $N$ term}. This is:
$$\on{term (4)} = g(R^g(x,y)x,y).$$
Note that because of the skew-symmetry of the Riemann tensor, this is a 
function of $x \wedge y$ also.

\noindent
\thetag{5} {\it The Ricci-curvature-like term}. The other curvature terms
are:
$$\tfrac12 (f_0^* g)^{ij}\! \Big(\!g(x,x)g(R^g(y,f_i)y,f_j)
-2g(x,y)g(R^g(x,f_i)y,f_j)+g(y,y)g(R^g(x,f_i)x,f_j)\! \Big)$$
If $V$ and $W$ are two perpendicular subspaces of the tangent space $TN_p$
at a point $p$, then we can define a `cross Ricci curvature' $\on{Ric}(V,W)$
in terms of bases $\{v_i\}, \{w_j\}$ of $V$ and $W$ by:
$$ \on{Ric}(V,W) = g^{ij} g^{kl} g(R^g(v_i,w_k)v_j,w_l).$$
Then this term factors as:
$$\on{term (5)} = \|x \wedge y\|_g^2 \on{Ric}(TM,\on{span}(x,y)).$$

\noindent
\thetag{6\,-7} {\it Terms involving the covariant derivative of $a$}. It is
remarkable that, so far, every term in the curvature tensor of $B_i$ vanishes if
$x \wedge y \equiv 0$, e.g., if the codimension of $N$ in $M$ is one! Now we have 
the terms:
\begin{align*} 
(f_0^*g)^{ij}\Big( &g(x,y)g(\nabla^\bot_i x,\nabla^\bot_j y) - \tfrac12
g(x,x)g(\nabla^\bot_i y, \nabla^\bot_j y) - 
\tfrac12 g(y,y)g(\nabla^\bot_i x, \nabla^\bot_j x) \\
&-g(x,\nabla^\bot_i x) g(y, \nabla^\bot_j y) 
-g(x,\nabla^\bot_i y) g(y, \nabla^\bot_j x)\\
&+g(x,\nabla^\bot_i y) g(x, \nabla^\bot_j y) 
+g(y,\nabla^\bot_i x) g(y, \nabla^\bot_j x)\Big).
\end{align*}
To understand this expression, we need a linear algebra computation, namely
that if $a,b,a',b' \in {\mathbb R}^n$, then:
\begin{multline*} 
<a,b><a',b'> - <a,a'><b,b'> -<a,b'><b,a'> - \\
-\tfrac12 <a,a><b',b'> -\tfrac12 <b,b><a',a'> + <a,b'>^2 + <b,a'>^2 =
\\
= \tfrac12(<a,b'>-<b,a'>)^2 - \tfrac12 \| a \wedge b' - b \wedge a'\|^2
\end{multline*}
Note that the term $g(x,\nabla^\bot y)$ (without an $i$) is a section of $\Om^1_M$
and the sum over $i$ and $j$ is just the norm in $\Om^1_M$, so the
above computation applies and the expression splits into 2 terms:
\begin{align*} 
\on{term (6)} &= -\tfrac12 \|(g(x,\nabla^\bot y)-g(y,\nabla^\bot x)\|^2_{\Om^1_M} 
\le 0\\
\on{term (7)} &= \tfrac12 \|x \wedge \nabla^\bot y - 
y \wedge \nabla^\bot x\|^2_{\Om^1_M\otimes \wedge^2 N(f)} \ge 0.
\end{align*}

Altogether, we get that the Riemann curvature of $B_i$ is the integral over
$M$ of the sum of the above 7 terms. We have the Corollary:

\begin{kor*} If the codimension of $M$ in $N$ is one, then all sectional
curvatures of $B_i$ are non-negative. For any codimension, sectional
curvature in the plane spanned by $x$ and $y$ is non-negative if $x$ and $y$ are
parallel, i.e., $x \wedge y = 0$ in $\bigwedge^2 T^*N$.
\end{kor*}

In general, the negative terms in the curvature tensor (giving positive 
sectional curvature) are clearly connected with the vanishing of geodesic 
distance: in some directions the space wraps up on itself in 
tighter and tighter ways. However, in codimension two or more with a flat
ambient space $N$ (so terms (4) and (5) vanish), there seem to exist conflicting 
tendencies making $B_i$ close up or open up: terms (1), (2) and (6) give positive
curvature, while terms (3) and (7) give negative curvature. It would be interesting
to explore the geometrical meaning of these, e.g., for manifolds of space curves.

\section{Vanishing geodesic distance on groups of
diffeomorphisms}\label{nmb:5}

\subsection{The $H^0$-metric on groups of diffeomorphisms}\label{nmb:5.1}
Let $(N,g)$ be a smooth connected Riemannian manifold, and let
$\on{Diff}_c(N)$ be the group of all diffeomorphisms with
compact support on $N$, and let $\on{Diff}_0(N)$ be the
subgroup of those which are diffeotopic in $\on{Diff}_c(N)$ to the identity;
this is the connected component of the identity in $\on{Diff}_c(N)$,
which is a regular Lie group in the sense of \cite{4}, section~38, see 
\cite{4}, section~42. The Lie algebra is $\X_c(N)$, the space of all 
smooth vector fields with compact support on $N$, with the negative of the
usual bracket of vector fields as Lie bracket. Moreover, $\on{Diff}_0(N)$
is a simple group (has no nontrivial normal subgroups), see \cite{E}, 
\cite{T}, \cite{Ma}. 
The {\it right invariant} $H^0$-metric 
on $\on{Diff}_0(N)$ is then given as follows, where $h,k:N\to TN$ are 
vector fields with compact support along $\ph$ and where 
$X=h\o\ph\i, Y=k\o\ph\i\in\X_c(N)$: 
\begin{align}
G^0_\ph(h,k) &= \int_N g(h,k)\on{vol}(\ph^*g)
= \int_N g(X\o\ph,Y\o\ph)\ph^*\on{vol}(g)
\notag\\&
= \int_N g(X,Y)\on{vol}(g)
\tag{1}\end{align}

\begin{thm}\label{nmb:5.2}
Geodesic distance on $\on{Diff}_0(N)$ with respect to the $H^0$-metric vanishes.
\end{thm}

\begin{demo}{Proof}
Let $[0,1]\ni t\mapsto\ph(t,\quad)$ be a smooth curve in
$\on{Diff}_0(N)$ between $\ph_0$ and $\ph_1$. Consider the
curve $u=\ph_t\o\ph\i$ in $\X_c(N)$, the right logarithmic
derivative. Then for the length and the energy we have:  
\begin{align}
L_{G^0}(\ph)&=\int_0^1  \sqrt{\int_N \|u\|^2_g\on{vol(g)}}\;dt
\tag{1}
\\
E_{G^0}(\ph)&=\int_0^1  \int_N \|u\|^2_g\on{vol(g)}\,dt
\tag{2}
\\
L_{G^0}(\ph)^2&\le E_{G^0}(\ph)
\tag{3} 
\end{align}

\noindent\thetag{4} Let us denote by $\on{Diff}_0(N)^{E=0}$ the set of all
diffeomorphisms $\ph\in\on{Diff}_0(N)$ with the following
property: For each $\ep>0$ there exists a smooth curve from the
identity to $\ph$ in $\on{Diff}_0(N)$ with energy $\le \ep$.

\noindent\thetag{5} {\it We claim that $\on{Diff}_0(N)^{E=0}$ coincides with
the set of all diffeomorphisms which can by reached from the
identity by a smooth curve of arbitraily short $G^0$-length.}
This follows by \thetag{3}.

\noindent\thetag{6} {\it We claim that $\on{Diff}_0(N)^{E=0}$
is a normal subgroup of $\on{Diff}_0(N)$.}
Let $\ph_1\in\on{Diff}_0(N)^{E=0}$ and $\ps\in\on{Diff}_0(N)$.
For any smooth curve $t\mapsto\ph(t,\quad)$ from the identity
to $\ph_1$ with energy $E_{G^0}(\ph)<\ep$ we have
\begin{align*}
E_{G^0}&(\ps\i\o\ph\o\ps) 
  = \int_0^1 \int_N\|T\ps\i\o\ph_t\o\ps\|_g^2
  \on{vol}((\ps\i\o\ph\o\ps)^*g)
\\&
\le \sup_{x\in N}\|T_x\ps\i\|^2 \cdot 
  \int_0^1 \int_N\|\ph_t\o\ps\|_g^2 (\ph\o\ps)^*\on{vol}((\ps\i)^*g)
\\&
\le \sup_{x\in N}\|T_x\ps\i\|^2 \cdot
  \sup_{x\in N}\frac{\on{vol}((\ps\i)^*g)}{\on{vol}(g)}\cdot
  \int_0^1 \int_N\|\ph_t\o\ps\|_g^2\,
  (\ph\o\ps)^*\on{vol}(g)
\\&
\le \sup_{x\in N}\|T_x\ps\i\|^2 \cdot
  \sup_{x\in N}\frac{\on{vol}((\ps\i)^*g)}{\on{vol}(g)}\cdot
  E_{G^0}(\ph).
\end{align*}
Since $\ps$ is a diffeomorphism with compact support, the two
suprema are bounded.
Thus $\ps\i\o\ph_1\o\ps\in\on{Diff}_0(N)^{E=0}$.
 
\noindent\thetag{7} {\it We claim that  
$\on{Diff}_0(N)^{E=0}$
is a non-trivial subgroup.} In view of the simplicity of
$\on{Diff}_0(N)$ mentioned in \ref{nmb:5.1} this concludes the
proof. 

It remains to find a non-trivial diffeomorphism in $\on{Diff}_0(N)^{E=0}$. 
The idea is to use compression waves. The basic case is this: take
any non-decreasing smooth function $f:{\mathbb R} \rightarrow {\mathbb R}$
such that $f(x)\equiv 0$ if $x \ll 0$ and $f(x) \equiv 1$ if $x \gg 0$.
Define 
$$\ph(t,x) = x + f(t-\la x)$$ 
where $\la < 1/\max(f')$. Note that 
$$\ph_x(t,x) = 1 - \la f'(t-\la x) > 0,$$ 
hence each map $\ph(t,\quad)$ is a diffeomorphism of $\mathbb R$ and 
we have a path in the group of diffeomorphisms of $\mathbb R$. These 
maps are not the identity outside a compact set however. In fact, 
$\ph(x)=x+1$ if $x \ll 0$ and $\ph(x)=x$ if $x \gg 0$. As 
$t \rightarrow -\infty$, the map $\ph(t,\quad)$ approaches
the identity, while as $t \rightarrow +\infty$, the map approaches
translation by 1. This path is a moving compression wave which pushes all
points forward by a distance 1 as it passes. We calculate its energy between
two times $t_0$ and $t_1$:
\begin{align*}
E_{t_0}^{t_1}(\ph) &= \int_{t_0}^{t_1} \int_{\mathbb R} 
\ph_t(t,\ph(t,\quad)^{-1}(x))^2 dx\, dt 
= \int_{t_0}^{t_1} \int_{\mathbb R} 
\ph_t(t,y)^2 \ph_y(t,y) dy\, dt  \\
&= \int_{t_0}^{t_1} \int_{\mathbb R} 
f'(z)^2 \cdot (1-\la f'(z)) dz/\la\, dt \\
& \le \frac{\max f'^2}{\la} \cdot (t_1-t_0) \cdot
\int_{\on{supp}(f')} (1-\la f'(z)) dz 
\end{align*}
If we let $\la = 1-\ep$ and consider the specific $f$ given by the convolution
$$f(z) = \max(0,\min(1,z))\star G_\ep(z),$$ 
where $G_\ep$ is a smoothing kernel supported on $[-\ep,+\ep]$, then the
integral is bounded by $3\ep$, hence
$$ E_{t_0}^{t_1}(\ph) \le (t_1-t_0)\tfrac{3\ep}{1-\ep}.$$

We next need to adapt this path so that it has compact support. To do this
we have to start and stop the compression wave, which we do by giving 
it variable length. Let:
$$ f_\ep(z,a) = \max(0,\min(a,z)) \star (G_\ep(z)G_\ep(a)).$$
The starting wave can be defined by:
$$ \ph_\ep(t,x) = x + f_\ep(t-\la x, g(x)),\quad \la <1, \quad 
g \text{ increasing}.$$
Note that the path of an individual particle $x$ hits the wave at $t=\la x-\ep$
and leaves it at $t=\la x + g(x) + \ep$, having moved forward to $x+g(x)$. 
Calculate the derivatives:
\begin{align*}
(f_\ep)_z &= I_{0 \le z \le a} \star (G_\ep(z)G_\ep(a)) \in [0,1]\\
(f_\ep)_a &= I_{0 \le a \le z} \star (G_\ep(z)G_\ep(a)) \in [0,1] \\
(\ph_\ep)_t &= (f_\ep)_z(t-\la x, g(x))\\
(\ph_\ep)_x &= 1 - \la (f_\ep)_z(t-\la x, g(x))
+(f_\ep)_a(t-\la x, g(x))\cdot g'(x) > 0.
\end{align*}
This gives us:
\begin{align*}
E_{t_0}^{t_1}(\ph) &= \int_{t_0}^{t_1} \int_{\mathbb R} (\ph_\ep)_t^2
(\ph_\ep)_x dx \, dt \\
&\le \int_{t_0}^{t_1} \int_{\mathbb R} (f_\ep)_z^2(t-\la x,g(x))\cdot
(1-\la (f_\ep)_z(t-\la x, g(x)))dx\, dt\\
&\quad+\int_{t_0}^{t_1} \int_{\mathbb R}(f_\ep)_z^2(t-\la x,g(x))\cdot
(f_\ep)_a(t-\la x, g(x)) g'(x) dx\, dt
\end{align*}
The first integral can be bounded as in the original discussion.
The second integral is also small because the support of the $z$-derivative
is $-\ep \le t-\la x \le g(x)+\ep$, while the support of the $a$-derivative
is $-\ep \le g(x) \le t-\la x + \ep$, so together $|g(x)-(t-\la x)| \le \ep$.
Now define $x_1$ and $x_2$ by $g(x_1)+\la x_1 = t+\ep$ and 
$g(x_0)+\la x_0 = t-\ep$.
Then the inner integral is bounded by
$$ \int_{|g(x)+\la x -t|\le \ep} g'(x) dx = g(x_1)-g(x_0) \le 2\ep,$$
and the whole second term is bounded by ${2\ep}(t_1-t_0)$. Thus the
length is $O(\ep)$.

The end of the wave can be handled by playing the beginning backwards. If
the distance that a point $x$ moves when the wave passes it is to be $g(x)$,
so that the final diffeomorphism is $x \mapsto x+g(x)$ then let $b = \max(g)$ 
and use the above definition of $\ph$ while $g' > 0$. The modification
when $g' < 0$ (but $g' > -1$ in order for $x \mapsto x+g(x)$ 
to have positive derivative) is given by:
$$ \ph_\ep(t,x) = x + f_\ep(t-\la x -(1-\la)(b-g(x)),g(x)).$$
A figure showing the trajectories $\ph_\ep(t,x)$ for sample values of $x$ is
shown in the figure above.

\begin{figure} 
\begin{center}
\epsfig{file=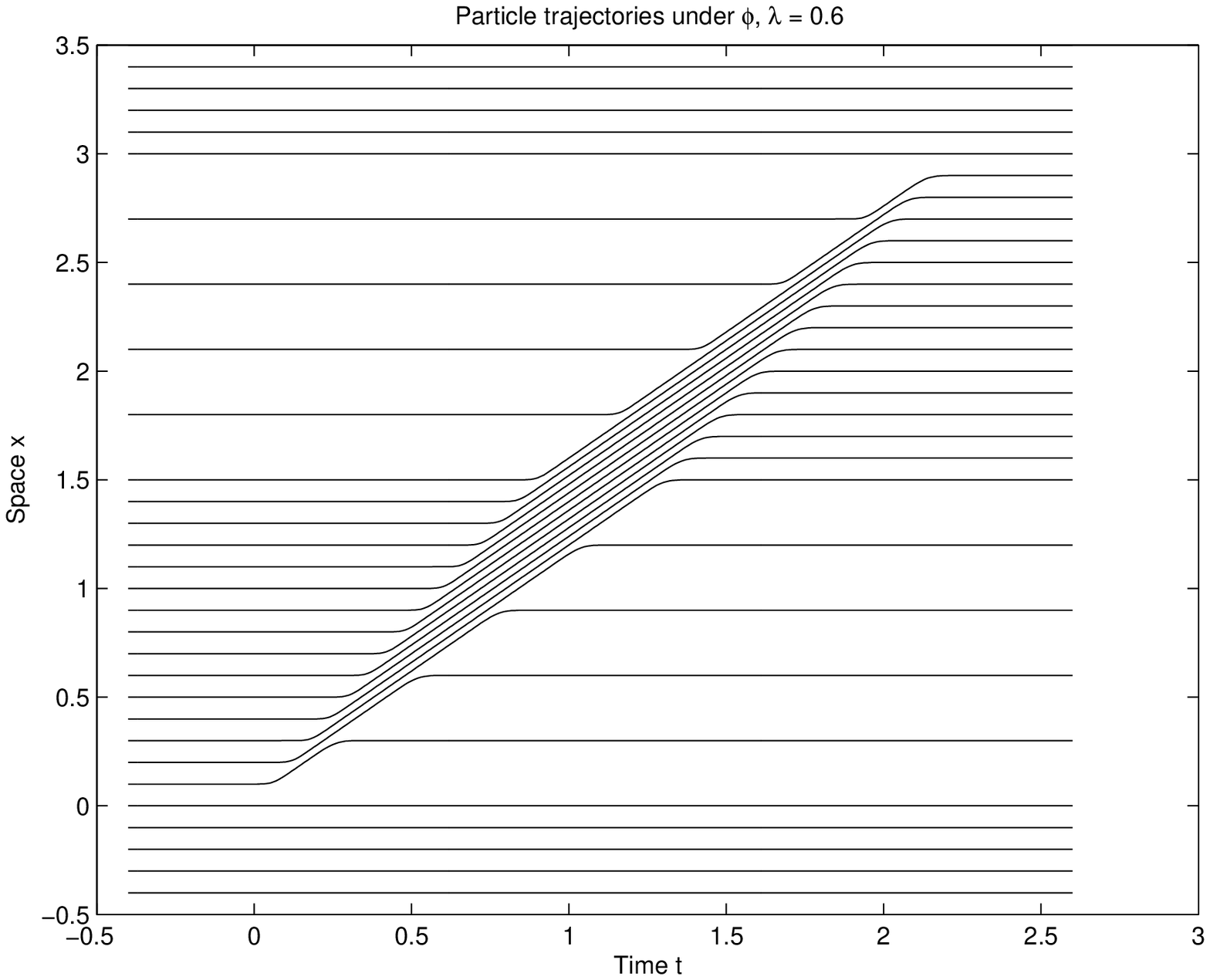,width=4in,height=3in}
\end{center}
\end{figure}

It remains to show that $\on{Diff}_0(N)^{E=0}$ is a nontrivial subgroup 
for an arbitrary Riemannian manifold. We choose a piece of a unit
speed geodesic containing no conjugate points in $N$ and Fermi coordinates 
along this geodesic; so we can assume that we are in an open set in 
$\mathbb R^m$ which is a tube around a piece of the $u^1$-axis. Now we 
use a small bump function in the the slice orthogonal to the $u^1$-axis 
and multiply it with the construction from above for the coordinate $u^1$. 
Then it follows that we get a nontrivial diffeomorphism in 
$\on{Diff}_0(N)^{E=0}$ again. 
\qed\end{demo}

\subsection*{Remark} Theorem \ref{nmb:5.2} can possibly be proved directly without
the help of the simplicity of $\on{Diff}_0(N)$. For $N=\mathbb R$ one can
use the method of \ref{nmb:5.2},~\thetag{7} in the parameter space of a
curve, and for general $N$ one can use a Morse function on $N$ to produce a
special coordinate for applying the same method, as we did in the proof of
theorem \ref{nmb:3.2}. 

\subsection{Geodesics and sectional curvature for $G^0$ on $\on{Diff}(N)$}
\label{nmb:5.3}
According to Arnold \cite{A}, see  \cite{MR},~3.3, for a right invariant
weak Riemannian metric $G$ on an (possibly infinite dimensional) Lie group
the geodesic equation and the curvature are given in terms of the adjoint operator 
(with respect to $G$, if it exists) of the Lie bracket
by the following formulas:
\begin{align*} 
u_t&=-\on{ad}(u)^*u, \quad u = \ph_t\o \ph\i
\\
G(\on{ad}(X)^*Y,Z)&:=G(Y,\on{ad}(X)Z)
\\
4G(R(X,Y)X,Y) &=  3G(\on{ad}(X)Y,\on{ad}(X)Y)                   
- 2G(\on{ad}(Y)^*X,\on{ad}(X)Y)  
\\&\quad
- 2G( \on{ad}(X)^*Y,\on{ad}(Y)X )                   
+ 4G( \on{ad}(X)^* X,\on{ad}(Y)^* Y )   
\\&\quad
- G( \on{ad}(X)^* Y+\on{ad}(Y)^* X,\on{ad}(X)^* Y+\on{ad}(Y)^* X ) 
\end{align*}
In our case, for $\on{Diff}_0(N)$, we have $\on{ad}(X)Y=-[X,Y]$ 
(the bracket on the Lie algebra $\X_c(N)$ of vector fields with
compact support is the negative of the usual one), and:
\begin{align*}
G^0(X,Y) &= \int_N g(X,Y)\on{vol}(g)
\\
G^0(\on{ad}(Y)^*X,Z)&= G^0(X,-[Y,Z]) = \int_N g(X,-\L_YZ)\on{vol}(g) 
\\&
= \int_N g\Bigl(\L_YX+(g\i\L_Yg)X+\on{div}^g(Y)X,Z\Bigl)\on{vol}(g)
\\
\on{ad}(Y)^*&= \L_Y+g\i\L_Y(g)+\on{div}^g(Y)\on{Id}_TN = \L_Y +\be(Y),
\end{align*}
where the tensor field $\be(Y)=g\i\L_Y(g)+\on{div}^g(Y)\on{Id}:TN\to TN$ 
is self adjoint with respect to $g$. 
Thus the geodesic equation is 
$$
u_t = -(g\i\L_u(g))(u)-\on{div}^g(u)u=-\be(u)u, \qquad u = \ph_t\o \ph\i.
$$
The main part of the sectional curvature is given by:
\begin{align*} 
4G(&R(X,Y)X,Y) =
\\&
= \int_N\Bigl( 3\|[X,Y]\|_g^2
+ 2g((\L_Y+\be(Y))X,[X,Y])  + 2g((\L_X+\be(X))Y,[Y,X] )                   
\\&\qquad\qquad
+ 4g(\be(X)X,\be(Y)Y) -\|\be(X)Y+\be(Y)X\|_g^2
\Bigr)\on{vol}(g)
\\&
= \int_N\Bigl(-\|\be(X)Y-\be(Y)X+[X,Y]\|_g^2
-4g([\be(X),\be(Y)]X,Y)
\Bigr)\on{vol}(g)
\end{align*}
So sectional curvature consists of a part which is visibly non-negative, and
another part which is difficult to decompose further. 

\subsection{Example: Burgers' equation}\label{nmb:5.4}
For $(N,g)=(\mathbb R,\on{can})$ or $(S^1,\on{can})$ the geodesic equation is 
{\it Burgers' equation} \cite{B}, a completely integrable infinite dimensional
system,
$$
u_t=-3u_x\,u, \qquad u = \ph_t\o \ph\i
$$
and we get $G^0(R(X,Y)X,Y)=-\int [X,Y]^2\,dx$ so that all sectional
curvatures are non-negative. 

\subsection{Example: $n$-dimensional analog of Burgers' equation}\label{nmb:5.5}
For $(N,g)=(\mathbb R^n,\on{can})$ or $((S^1)^n,\on{can})$ we have:
\begin{align*}
(\on{ad}(X)Y)^k&=\sum_i((\p_iX^k)Y^i-X^i(\p_iY^k))
\\
(\on{ad}(X)^* Z)^k &= \sum_{i} \Bigl( 
     (\p_kX^i)Z^i + (\p_iX^i)Z^k + X^i(\p_iZ^k)\Bigr), 
\end{align*}
so that the geodesic equation is given by 
$$ 
\p_t u^k = - (\on{ad}(u)^\top u)^k  
=  -\sum_{i} \Bigl((\p_ku^i)u^i + (\p_iu^i)u^k + u^i(\p_iu^k)\Bigr), 
$$
called the basic Euler-Poincar\'e equation (EPDiff) in \cite{HRTY},
the $n$-dimensional analog of Burgers' equation.

\subsection{Stronger metrics on $\on{Diff}_0(N)$}\label{nmb:5.6}
A very small strengthening of the weak Riemannian $H^0$-metric on 
$\on{Diff}_0(N)$ makes it into a true metric. We define the stronger 
right invariant semi-Riemannian metric by the formula:
$$G^A_\ph(h,k)= \int_N (g(X,Y)+A\on{div}_g(X).\on{div}_g(Y)) \on{vol}(g).$$
Then the following holds:

\begin{thm}\label{nmb:5.7}
For any distinct diffeomorphisms $\ph_0, \ph_1$, the infimum of the
lengths of all paths from $\ph_0$ to $\ph_1$ with respect to $G^A$
is positive.
\end{thm}

This implies that the metric $G^0$ induces positive geodesic
distance on the subgroup of volume preserving diffeomorphism since it
coincides there with the metric $G^A$.

\begin{demo}{Proof}
Let $\ps_1 = \ph_0 \circ \ph_1^{-1}$. If $\ph_0 \ne \ph_1$, there are two
functions $\rh$ and $f$ on $N$ with compact support such that:
$$ \int_N \rh(y) f(\ps_1(y)) \on{vol}(g)(y) \ne 
\int_N \rh(y) f(y) \on{vol}(g)(y). $$
Now consider any path $\ph(t,y)$ between the two maps with derivative 
$u = \ph_t \circ \ph^{-1}$. Inverting the diffeomorphisms (or switching
from a Lagrangian to an Eulerian point of view), let 
$\ps(t,\quad) = \ph(0,\quad)\circ \ph(t,\quad)^{-1}$. Then $\ps_t = -\D\ps(u)$
and we have:
\begin{align*}
\int_N &\rh(y) f(\ps_1(y)) \on{vol}(g)(y) -\int_N \rh(y) f(y)\on{vol}(g)(y) =
\\&
=\int_0^1 \int_N \rh(y) \p t f(\ps(t,y) \on{vol}(g)(y) dt
= \int_0^1\int_N \rh(y) (df \circ \ps)(\ps_t(t,y))\on{vol}(g)(y)\, dt
\\&
= \int_0^1\int_N \rh(y) (\D f \circ \ps)(-\D\ps(u(t,y)))\on{vol}(g)(y) dt
\end{align*}
But $\on{div}((f\circ \ps) \cdot \rh u)=(f\circ \ps) \cdot \on{div}(\rh u)
+ (\D f\circ \ps)( \D\ps(\rh u))$. The integral of the left hand side is 0, hence:
\begin{align*}
\Big|\int_N \rh(y) f(\ps_1(y)) \on{vol}(g)(y) &- \int_N \rh(y) f(y)
\on{vol}(g)(y)\Big|
\\& 
= \Big|\int_0^1 \int_N (f\circ \ps)\on{div}(\rh u) \on{vol}(g) dt\Big|
\\&
\le \sup(|f|)\int_0^1 \sqrt{\int_N C_\rh \|u\|^2 + 
C'_\rh|\on{div}(u)|^2 \on{vol}(g)}dt
\end{align*}
for constants $C_\rh, C'_\rh$ depending only on $\rh$. Clearly the right hand 
side is a lower bound for the length of any path from $\ph_0$ to $\ph_1$.
\qed\end{demo}

\subsection{Geodesics for $G^A$ on
$\on{Diff}(\mathbb R)$}\label{nmb:5.8}
See \cite{CH} and \cite{Mis}. 
We consider the groups $\on{Diff}_c(\mathbb R)$ or $\on{Diff}(S^1)$ with Lie  
algebras $\X_c(\mathbb R)$ or $\X(S^1)$ with Lie bracket  
$\on{ad}(X)Y=-[X,Y]=X'Y-XY'$.
The $G^A$-metric equals the $H^1$-metric on $\X_c(\mathbb R)$, and we have: 
\begin{align*} 
G^A(X,Y) &= \int_{\mathbb R} (XY + AX'Y')dx=
\int_{\mathbb R}X(1-A\p_x^2)Y\,dx,\\ 
G^A(\on{ad}(X)^* Y&, Z)= \int_{\mathbb R} (YX'Z-YXZ'+AY'(X'Z-XZ')') dx\\ 
&= \int_{\mathbb R}  
Z(1-\p_x^2)(1-\p_x^2)\i(2YX'+Y'X-2AY''X'-AY'''X) dx,\\ 
\on{ad}(X)^* Y &= (1-\p_x^2)\i(2YX'+Y'X-2AY''X'-AY'''X)\\ 
\on{ad}(X)^* &= (1-\p_x^2)\i(2X'+X\p_x)(1-A\p_x^2)
\end{align*}
so that the geodesic equation in 
Eulerian representation   
$u=(\p_t f)\o f\i \in \X_c(\mathbb R)$ or $\X(S^1)$ is
\begin{align*} 
\p_t u &= - \on{ad}(u)^* u
=  -(1-\p_x^2)\i(3uu'-2Au''u'-Au'''u),\text{ or }\\ 
u_t-u_{txx}&= Au_{xxx}.u+2Au_{xx}.u_x-3u_x.u,
\end{align*}
which for $A=1$ is the {\it Camassa-Holm equation} \cite{CH}, another 
completely integrable infinite dimensional Hamiltonian system. 
Note that here geodesic distance is a well defined metric describing the
topology.

\end{document}